\documentclass[12pt,twoside,a4paper]{amsart}
\usepackage{amssymb,amsthm,amsmath,xy,lscape}
\usepackage{a4wide}
\xyoption{all}
\usepackage[swedish,english]{babel}
\usepackage[latin1]{inputenc}
\usepackage[hypertex,
colorlinks=true, urlcolor=blue, pagecolor=black,
linkcolor=black, citecolor=black, filecolor=black, menucolor=black, bookmarks=true, bookmarksnumbered=true,
pdfpagelayout=TwoColumnRight, pdfstartview=FitH, plainpages=false,
pdfpagelabels]{hyperref}
\usepackage{multicol}

\font\eightsy=cmsy8
\font\tensy=cmsy10
\font\eightib=cmmib8

\theoremstyle{plain}
\newtheorem{theorem}{Theorem}[section]
\newtheorem{proposition}[theorem]{Proposition}

\newtheorem*{varprop}{Proposition}

\theoremstyle{definition}
\newtheorem{definition}[theorem]{Definition}
\newtheorem*{vardef}{Definition}

\def\DEFhatH{{
\hat{\hbox{\bf H}}}}
\def\hatH{$\DEFhatH$}
\def\mhatH{\DEFhatH}

\begin{document}

\title{Algebraic Topology 2.0}
\author{G\"oran Fors}
\address{\noindent\hskip-0.5cm Department of Mathematics, University of Stockholm, 
Stockholm, Sweden. goranf@math.su.se}

\begin{abstract}{
The efficiency of contemporary algebraic topology is not optimal since the category of topological spaces can be made more algebraic by introducing a profoundly new $(-1)$-dimensional topological space \{$\wp$\} as a topological join unit. Thereby synchronizing the category of topological spaces with the structures within the contemporary category of simplicial complexes as well as with the structures within the algebraic categories.

In the category of topological spaces, the empty space $\emptyset$ has since long been given the role as a join unit - ad-hoc though. Since it is $\{\emptyset\}$, not $\emptyset$, that is the join unit within the category of simplicial complexes, that role of $\emptyset$ within general topology  has to be rectified.

This article presents an algebraization of Hausdorff's century old definition of the category of topological spaces as well as some useful algebraic topological consequences thereof.}
\end{abstract}

\maketitle
\tableofcontents

\def\myl{{\vcenter{\hrule width 0.13 true in height 0.005 true in}}}
\def\myline{$\myl$}
\def\mmyline{\myl}

\def\mysqcuplus{{\rlap{$\sqcup$}{\raise1.5pt\hbox{\sevenbf +}}}}
\def\mysqcupplus{\mysqcuplus}
\def\sqcupP{\ \!{{\rlap{$_{_{\!}}\sqcup
$}{\vbox{\moveleft2.8pt\hbox{{\ }\raise1.0pt\hbox{\fivebf P}}}}} } }
\def\mysqcupP{\sqcupP}
\def\Psqcup#1{\ \!{{\rlap{$_{_{\!}}\sqcup
$}{\vbox{\moveleft2.0pt\hbox{{\ }\raise1.0pt\hbox{\fivebf
\hbox{#1}}}}}} } }
\def\myPsqcup#1{\Psqcup#1~}

\def\mmysubsetneqq{\subsetneqq}
\def\mmysupsetneqq{\supsetneqq}

\def\nabl#1
{\font\mysyfont=cmsy#1
{
{\!_{\!}}^{_{_{\hbox{\mysyfont \char"72}}}}\!}}
\def\mynabla#1{$\nabl#1 $}
\def\mmynabla#1{\nabl#1 }

\def\barnabl#1{\!\bar{\ \!\nabl#1 }}
\def\mybarnabla#1{\barnabl#1 }

\def\thickl{{\vcenter{\hrule width 0.1 true in height 0.015 true in }}}
\def\thickline{$\thickl$}
\def\mthickline{\thickl}

\def\DEFhatH{{
\hat{\hbox{\bf H}}}}
\def\hatH{$\DEFhatH$}
\def\mhatH{\DEFhatH}

\def\eDEFhatH{{
\hat{\hbox{\eightbf H}}}}
\def\ehatH{$\eDEFhatH$}
\def\mehatH{\eDEFhatH}

\def\DEFsmallhatH{{
\hat{\hbox{\eightbf H}}}}
\def\smallhatH{$\DEFsmallhatH$}
\def\msmallhatH{\DEFsmallhatH}

\def\DEFminihatH{{
\hat{\hbox{\sixbf H}}}}
\def\minihatH{$\DEFminihatH$}
\def\mminihatH{\DEFminihatH}

\def\FFrame#1#2#3{                   %
     \vbox{\hrule height#2           
          \hbox{\vrule width#2       %
                 \hskip#1            
                 \vbox{\vskip#1{}    %
                       #3            
                       \vskip#1}     %
                 \hskip#1
                 \vrule width#2}
           \hrule height#2}}

\def\DEFhatbar{\hbox{$\mathbf{\hat{\vert}}$}}
\def\hatbar{$\DEFhatbar$}
\def\mhatbar{\DEFhatbar}

\def\DEFhatSigma{\hbox{$\mathbf{\hat{\Sigma}}$}}
\def\hatSigma{$\DEFhatSigma$}
\def\mhatSigma{\DEFhatSigma}

\def\DEFHatAst{{
{
{\hbox{$^{{\land}}$}}{\vbox{\hbox{$\hskip-5.5pt\ast$}}}}
}}
\def\hatast{$\DEFHatAst$}
\def\mhatast{\DEFHatAst}

\def\DEFsmallHatAst{{\ \!\!{\rlap{{\lower3.5pt\hbox
{\vbox{\moveright-0.5pt\hbox{$^{^{{\hbox{\fivesy \char"5E}}}}$}}}}}
{{{\hbox{\fivesy \char"03}}}}}\ \! }}
\def\topast{$\DEFsmallHatAst$}
\def\mtopast{\DEFsmallHatAst}

\def\DEFsmallHattAst{{\ \!\!{\rlap{{\lower2.0pt\hbox
{\vbox{\moveright-0.5pt\hbox{$^{^{{\hbox{\fivesy \char"5E}}}}$}}}}}
{{{\hbox{\sevensy \char"03}}}}}\ \! }}
\def\toppast{$\DEFsmallHattAst$}
\def\mtoppast{\DEFsmallHattAst}

\def\DEFCircast{{{\rlap{\lower4.5pt\hbox{$^{^{_{{_{\!}}\bigcirc}}}$}}
{\vbox{\moveright-0.5pt\hbox{$\ast$}}}}}}
\def\circast{$\DEFCircast$}
\def\mcircast{\DEFCircast}

\def\DEFsimpline{\bullet\!\!\!-\!-\!-\!\!\!\bullet}
\def\simplline{$\DEFsimpline$}
\def\msimplline{\DEFsimpline}

\def\DEFsimpsquare{{\hbox{\hsize1.23cm\hskip0.1cm
\raise0.7pt\vbox{\hrule width 1.0pt height 0.3cm depth
0.32pt}\hskip-0.48cm
\vbox{$_{\hbox{$\bullet\!\raise2.8pt\vbox{\hrule width 0.25cm height 0.7pt depth 0.32pt}\!\bullet$}}
^{^{\hbox{$\bullet\!\raise2.8pt\vbox{\hrule width 0.25cm height
0.7pt depth 0.32pt}\!\bullet$}}}$}
\hskip-0.5cm\raise0.7pt\vbox{\hrule width 1.0pt height 0.3cm depth
0.32pt}\hskip0.1cm
}
}}
\def\simplsquare{$\DEFsimpsquare$}
\def\msimplsquare{\DEFsimpsquare}

\def\DEFringhom#1{
{\lower4.0pt\hbox{$^{_{\hskip0.05cm\hbox{{\sevenbf #1}}
\atop{\hbox{$\cong$}}}}$}}}
\def\ringHom#1{$\DEFringhom#1$}
\def\mringHom#1{\DEFringhom#1}

\def\DEFringhombb#1{
{\hbox{$^{_{{\mathbb{#1}}
\atop{\hbox{$\cong$}}}}$}}}
\def\ringHombb#1{$\DEFringhombb#1$}
\def\mringHombb#1{\DEFringhombb#1}

\def\DEFsmallRinghom#1{
{\lower4.0pt\hbox{$^{_{\hbox{{\fivebf #1}}
\atop{{\raise1.4pt\hbox{{\eightbf =}}}
\hskip-0.25cm{\raise3.8pt\hbox{{\eightsy {\char"18}}}} }}}$}}}
\def\smallRingHom#1{$\DEFsmallRinghom#1$}
\def\msmallRingHom#1{\DEFsmallRinghom#1}

\def\DEFbfdel#1{\hbox{\tenbf#1$\hskip-3.0pt$#1$\hskip-3.0pt$#1}}
\def\bfDel#1{$\DEFbfdel#1$}
\def\mbfDel#1{\DEFbfdel#1}

\def\DEFbigbfdel#1{\hbox{\tenbf\big#1$\hskip-4.5pt$\big#1$\hskip-4.5pt$\big#1}}
\def\BfDel#1{$\DEFbigbfdel#1$}
\def\mBfDel#1{\DEFbigbfdel#1}

\def\DEFbfbrace#1{\hbox{$#1$\hskip-5.8pt$#1$\hskip-5.8pt$#1$}}
\def\bfBrace#1{$\DEFbfbrace#1$}
\def\mbfBrace#1{\DEFbfbrace#1}

\def\DEFbigbfbrace#1{\hbox{$\big#1$\hskip-6.75pt$\big#1$\hskip-6.75pt$\big#1$}}
\def\BfBrace#1{$\DEFbigbfbrace#1$}
\def\mBfBrace#1{\DEFbigbfbrace#1}

\def\DEFbfCupCap#1#2{\hbox{$#1$\hskip-0.33cm$\lower3pt\hbox{$_{_{#2}}$}$}}
\def\bfcupcap#1{$\DEFbfCupCap#1$}
\def\mbfcupcap#1{\DEFbfCupCap#1}

\def\mytinynabla{\lower1.0pt\hbox{$^{_{_{\hbox{\sevensy \char"72} }}}$}}

\def\LHS{{\bf L\hatH S}\hskip2pt}

\section{Introduction}\label{SecP0:1}
We give four quotations to catch some aspects of the foundation of {\it Algebraic topology}.

\smallskip
First, \cite{Ste} N. E. Steenrod: {\it A Convenient Category of Topological Spaces} (1967) pp. 1-2:
\begin{quote}
(p.~1) For many years, algebraic topologists have been laboring under the handicap of not knowing in which category of spaces they should work. ...\\
(p.~2) In this paper, we propose as convenient the category of spaces we shall call compactly generated.  Such a space is a Hausdorff space with the property that each subset that intersects every compact set in a closed set is itself a closed set.
\end{quote}

\smallskip
When J. H. Ewing reviewed William S. Massey's monograph
{\it Homology and cohomology theory} in \cite{Ew} (1979) 
he made quite a few reflections - some of which follow below:
\begin{quote}
Algebraic topology attempts to solve topological problems using algebra. To do so requires some sort of machine which produces the ``algebraic image" of topology, and it is the machine itself on which topologists often spend most of their time, first carefully building and then diligently refining.
Historically, the first such machine was ordinary homology and cohomology theory. (The word ``ordinary" is not a slur~it means homology and cohomology defined from an algebraic chain complex as opposed to "extraordinary" theories such as {\it K}-theory.)\\
\vdots\vskip-0.4cm
\dots I have stated that homology and cohomology are at the very foundation of algebraic topology, and, moreover, that topologists have been confused for some time about which
particular variant to use. Some readers are apt to conclude that,\\
(i) algebraic topology is a subject in great disarray,\\
(ii) algebraic topologists have certainly
been a careless and giddy lot to allow such a situation to continue, and\\
(iii) the news of this breakthrough%
({\tensy{\char"6E}{\char"6E}}
of Massey using results of Nöbeling to promote the use of the
{\it Steenrod  homology}{\tensy{\char"6E}{\char"6E}})%
will surely excite all topologists.\vskip0.0cm
Lest the reader be led astray it should be forcefully pointed out that most of the time most topologists only consider spaces for which all theories do agree. Many
machines of algebraic topology breakdown when applied to general spaces;
\vskip-0.2cm
\dots
\vskip-0.01cm
Yet, as we indicated above, algebraic topology has developed far beyond its roots and in the process much debris has been left behind. There is a great need for simplification, unification and illuminating exposition; it is a need which is
largely unmet. Perhaps it's time for some to stop partying just long enough to tidy up.
\end{quote}

F. Quinn lay out the development of the subject {\it controlled topology} in \cite{Q} (2002) pp.~1-2:
\begin{quote}
\noindent
{\bf 1.1 Locating the subject.} In the first half of the 20th century topology had two main branches: point-set topology, concerned with local properties (separation, connectedness, dimension theory etc); and algebraic topology, concerned with definition and detection of global structure (homology, characteristic classes, etc.).
In the 50s and 60s the algebraic branch split into homotopy theory and geometric topology....
Controlled topology began in the late 1970s and 80s as a way to apply the constructive techniques of geometric topology to local questions more typical of point-set topology....
\end{quote}
 The ideational roots of "category theory" is indicated in
\cite{McL} (2005) - C. McLarty: {\it Saunders MacLane (1909-2005):
His Mathematical Life and Philosophical Works}. 
In p.~238, after mentioning some connections and similarities between Eliakim Hastings Moore (1862-1932) and Saunders MacLane, McLarty points out that Moore had a principle stated as follows:

\vskip-0.17cm\vskip-0.17cm
\begin{quote}
\noindent
-``The existence of analogies between central features of
various theories implies the existence of a general abstract
theory which underlies the particular theories and unifies
them with respect to those central features". (Moore [1908], p. 98)\\
Category theory was created to find the unity behind deep
analogies between topology and algebra. 
Today it is a standard format for giving such `general abstract' theories.
\end{quote}

Learning mathematics the Bourbaki way, through {\it structures and axioms},
has its advantages as stated in 
\cite{McC} (2005)
\htmladdnormallink{{\it Bourbaki and Algebraic Topology} }{http://mathserver.vassar.edu/faculty/McCleary/BourbakiAlgTop.pdf}
by J. McCleary. Then again, modern mathematics is woven around {\it category theory}.
When combinatorialists in the 70th begun to use the Stanley-Reisner ring-construction, def.~\ref{defStR-ring} p.~\pageref{defStR-ring}, to solve combinatorial problems, it was immediately clear that the \underbar{classical} category of simplicial complexes was inadequate for their purposes and that it had to be - and was - rectified to allow a vertex-free $(-1)$-dimensional simplicial join-unit $\{\emptyset\}$ matching the unit object ``${\bf k}\circ{\rm Id}$" (written ${\bf k}$ for short) w.r.t. the ${\bf k}$-tensor product among ${\bf k}$-algebras while the empty simplicial complex $\emptyset$ now became to match the trivial ring {\bf 0}, cf. Examples {\bf ii} and {\bf iii} p.~\pageref{Ex1:03}.

Now, quite obviously, the \underbar{classical} topological category that was underlying the above mentioned creation of ``category theory" is equally inadequate to match the algebra-structure as was the category of \underbar{classical} simplicial complexes. Topology need two \underbar{point-free} spaces - the topological join-zero $\emptyset$ and a new join-unit \{$\wp$\} to match even the contemporary combinatorial structure, with its two \underbar{vertex-free} complexes - simplicial join-zero $\emptyset$ and the simplicial join-unit $\{\emptyset\}$ invoked under the mentioned algebraization in the 1970th.

General topology were not equipped with any second point-free structure-vital object \{$\wp$\} until the beginning of the 1990th when the early draft of this author's thesis addressed that issue and provided that last piece of the above mentioned synchronization-puzzle regarding general topology, combinatorics and algebra.
Now, the {\it realisation} of augmental simplicial complexes can be done through a \underbar{faithful} covariant functor as in the purely classical setting.

The switch from {\it Hausdorff's classical category of topological spaces and continuous functions} to the {\it augmental category} is further motivated by the following examples, starting with a quotation from \cite{Wh2} G.W. Whitehead: \htmladdnormallink{{\it Generalized Homology Theories}
}{http://www.ams.org/journals/tran/1962-102-02/S0002-9947-1962-0137117-6/S0002-9947-1962-0137117-6.pdf}(1962) p.~228:

\vskip-0.1cm\vskip-0.1cm
\begin{quote}

\noindent
2. {\bf Preliminaries.} Let ${\mathcal W}_0$ be the category of spaces with base-point having the homotopy type of a CW-complex. More precisely, an object of ${\mathcal W}_0$ is a
space $X$ with base-point $x_0$, such that there exists a CW-complex $K$ with base-point $k_0$ and a homotopy equivalence of the pairs $(X, \{x_0\})$ and $(K, \{k_0\})$; and a map of ${\mathcal W}_0$ is a continuous, base-point preserving map.

Let ${\mathcal W}$ be the category of spaces (without distinguished base-point) having
the homotopy type of a CW-complex. Let $P$ be a fixed space consisting of exactly one point $p_o$. If $X\in {\mathcal W}$, let $X^+$ be the topological sum of $X$ and $P$; then ($X^+, p_0)$ is an object of ${\mathcal W}_0$.
If $X, Y\in {\mathcal W}$ and $f: X\rightarrow Y$, then $f$ has a unique
extension $f^+: X^+\rightarrow Y^+$ such that $f^+(p_0) =p_0$, and $f^+$ is a map in ${\mathcal W}_0$. The correspondences
$X\rightarrow X^+$, $f\rightarrow f^+$ define a functor $^+: {\mathcal W}\rightarrow {\mathcal W}_0$. Evidently we may
regard ${\mathcal W}$ as a subcategory of ${\mathcal W}_0$.
\end{quote}

The last quote describes a familiar routine from Homotopy Theory providing all free spaces with a common base point - essentially identical to what is done in \cite{Wh1} G.W. Whitehead: {\it Elements of Homotopy Theory}, GTM~61, Springer (1978) p.~103.

When passing from the classical category of topological spaces to the augmental one, the process is similar to the steps taken by Whitehead above except that one uses the now available \underbar{unique} $(-1)$-dimensional space $\{\wp\}$ instead of $ P=\{p_0\}$, while at the same time avoiding the choice of a {\it particular point} $p_0\in P=\{p_0\}\subset X$. N.B. $\{p_0\}$ is $0$-dimensional.

\cite{Wh3}G.W. Whitehead:%
\htmladdnormallink{\emph{Homotopy Groups of Joins and Unions}}{http://www.ams.org/journals/tran/1956-083-01/S0002-9947-1956-0080918-0/S0002-9947-1956-0080918-0.pdf} (1956) shows a related issue.
\begin{quote}
\noindent
(p~56.) The join of X with the empty set $\emptyset$ is $X$. (sic. An ad-hoc-convention invented by the general topologists.)\\
(p~57.) We consider the empty set as an ordered Euclidean
$(-1)$-simplex. With the usual definition of equivalence of singular simplexes, we define $\tilde{S}_p(X)$ to be the free abelian group generated by the singular $(p-1)$-simplexes in $X$ $(p = 0, 1, 2,\cdots)$ and define $\tilde{S}(X) = \sum_{p=0}^{\infty}\tilde{S}_p(X)$, with the usual boundary operator (the boundary of each singular $0$-simplex is the unique singular $(-1)$-simplex), to be the {\it augmental total singular complex} of X. If $A\subset X$, then $\tilde{S}(A)$ is a sub\-complex of $\tilde{S}(X)$ and we define $\tilde{H}_p(X, A)$ to be the $(p + 1)$st homology group of the complex $\tilde{S}(X)/\tilde{S}(A)$. Of course, $\tilde{H}_p(X,A)
= H_p(X, A)\ (p\geqq0)$, except that $\tilde{H}_0(X,\emptyset)$ is the {\it reduced} $0$-dimensional homology group of $X$.
\end{quote}

The claim above ``We consider the empty set as an ordered Euclidean $(-1)$-simplex" implies that the singular chain-complex for $\emptyset$ will have ${\rm Id}_\emptyset$ as a $(-1)$-dimensional generator, which contradicts the last sentence ``... that $\tilde{H}_0(X,\emptyset)$ is the reduced $0$-dimensional homology group of $X$", since then  $\tilde{S}_i(\emptyset)$ necessary equals the trivial group ${\bf 0}$ for every $i$.

Finally, we quote \cite{MacL} S. Mac~Lane: {\it Categories for the Working Mathematician} (1998);
\begin{quote}
\noindent
{\rm(p.175)}. {\bf The Simplicial Category:}
This category
{\eightib \char"01}
has as objects all {\it finite} {\it ordinal} {\it numbers} 
\hbox{\bf n}:=\hbox{\eightsy \char"66}0,\dots, n-1\hbox{\eightsy\char"67}
and as arrows
\hbox{\bf f}:\hbox{\bf n}\hbox{\eightsy \char"21}\hbox{\bf
n}$^{\prime}$ all (weakly) monotone functions: \dots\\
{\bf \noindent
\noindent {\rm (p.~178)}:}
The category $\hbox{\eightib\char"01}$ has a direct geometric interpretation by affine simplices, which give a functor $\Delta:$\hbox{\eightib\char"01}$\rightarrow {\bf Top}$ representing $\hbox{\eightib\char"01}$ as a subcategory of ${\bf Top}$. On objects $n$ of $\hbox{\eightib\char"01}$ take $\Delta_0$ to be the empty topological space, and $\Delta_{n+1}$ to be the ``standard" $n$-dimensional affine simplex - the subspace of Euclidean $\mathbf{\mathbb{R}}^{n+1}$ consisting of the following points
 $\Delta_n=\{p=(t_0,...,t_n\vert t_0\geqq0,...,t_n\geqq0),\sum{t_i}=1\};$\\
\end{quote}
\vskip-0.2cm\vskip-0.2cm
-``On objects $n$ of $\hbox{\eightib\char"01}$ take $\Delta_0$ to be the empty topological space,\dots" indicates that the empty space $\emptyset$ is to be regarded as a $(-1)$-dimensional topological space. Here, as well as in the three mentioned text-examples from  G.W. Whitehead's texts above, the author tempers with the dimensions of objects that today have well-established dimensions within contemporary combinatorics, i.e. $\dim\emptyset = -\infty$, $\dim\{\emptyset\} = -1$ and $\dim\{{\rm vertex}\} = 0$, where ``${\rm vertex}$" turns into ``${\rm point}$" in the category of topological spaces. A topological counterpart $\{\wp\}$ of $\{\emptyset\}$ has yet to be introduced within general topology, as is done below.
The category of $\{\wp\}$-extended topological spaces $\mathcal{C}_o$ is called  \hbox{\it the augmental topological category}.

The definitions of {\it continuous functions, standard simplexes, setminus, quotients, homology manifolds} etc within the augmental category $\mathcal{C}_o$, are properly reformulated - respecting and embedding those of the related concepts within the classical topological category $\mathcal{C}$.

The definitions of manifolds and their boundaries are given in p.~\pageref{quasim} and the following boundary formula for the join of two disjoint manifolds is established in Theorem~\ref{Th01P22} p.~\pageref{Th01P22}:
$${\rm Bd} ({M}_{1}\ast
{M}_{2})
=(({\rm Bd} {M}_{1})\ast {M}_{2})\cup ({M}_{1}\ast ({{\rm Bd}} {M}_{2})).
$$
This formula does not live classically. With ${\mathcal M} \ast \bullet$ being ``the cone of a Möbius band" and using the prime field $\mathbb{Z}_{\bf 3}$ with three elements as coefficient for the homology groups, the above formula gives (for quasi-manifolds);
${\rm Bd}_{_{{\mathbb{Z}_{\bf 3}}}}\!\!({\mathcal M} \ast \bullet)=$
$(({\rm Bd}\!_{_{{\mathbb{Z}_{\bf 3}}}}\!\!
{\mathcal M})\ast\bullet) \cup ({\mathcal M} \ast
{\rm Bd}_{_{{\mathbb{Z}_{\bf 3}}}}\!\!
\bullet)
=
(({\mathbb{S}}^1)\ast\bullet)\cup ({\mathcal M}
\ast
{\rm Bd}_{_{{\mathbb{Z}_{\bf 3}}}}\!\! \bullet)
=
\mathbb{B}^2 \cup
({\mathcal M}\ast \{\emptyset_o\})
=\mathbb{B}^2 \cup {\mathcal M}
= {{\mathbb{RP}}}^2\!\!:=$
\htmladdnormallink{{\it real projective plane}}
{http://en.wikipedia.org/wiki/Real_projective_plane}, where
$\mathbb{B}^n$ is the $n$-dimensional disk or
$n$-ball and $\mathbb{S}^n$ is the $n$-sphere.
$\emptyset, \{\emptyset_o\},$ and $0$-dimensional complexes
with either one, ${\bullet},$ or two, ${\bullet}{\bullet},$ vertices are the only manifolds in dimensions $\le0$.

${\rm Bd} \bullet=$``the boundary of the $0$-ball $\bullet$" $=\{\emptyset_o\} =\mathbb{S}^{-1}$. The double of $\bullet$ is the $0$-sphere $\bullet\bullet$. Both the $(-1)$-sphere
$\{\emptyset_o\}$ and the $0$-sphere $\bullet\bullet$ has, as preferred, empty boundaries.
{\bf Note} that
classically: ${\rm Bd} \bullet= {\rm Bd} \bullet\bullet = \emptyset$.
Now one is bound to conclude that the wrong category of topological spaces has been used within algebraic topology - all along.

The use of the category of augmental topological spaces and continuous functions instead of Hausdorff's classical category do carry the sign of a paradigm shift, since none of the following useful example-formulas within $\mathcal{C}_o$ from this article can live in the classical setting.

\smallskip
\noindent
{\bf Theorem 3.4.}\label{Th3.4} (See p.~\pageref{Th04}) {\it 
{\rm[The \htmladdnormallink{{\it Künneth formula} }
{http://en.wikipedia.org/wiki/K\%C3\%BCnneth_theorem} for Topological Joins; cp. \cite{Sp} p.~235.]}\\
If $\{X_{1}\ast Y_{2}, X_{2}\ \ast Y_{1}\}$ is an excisive couple in $X_{1}\ast Y_{1},$ ${\mathbf{R}}$ a ${\mathbf{PID}}$, $\mathbf{G} \mbox{ and } \mathbf{G}^\prime\
\mathbf{R}$-modules and $\mathrm{Tor}_1^{\mathbf{R}}(\mathbf{G},\mathbf{G}^{\prime})=\mbox{\bf0}.$
Then the functorial sequences below are $($non-naturally$)$ split exact};
$\begin{array}{lll}
{{\bf0}\longrightarrow \bigoplus_{i+j=q} \ \ [{\hat{\mathbf{H}}}_{_{i}}
           (X\!_{_1},X\!_{_2};{\bf G})
\otimes\!\!_{_{\mathbf{R}}}
        {\hat{\mathbf{H}}}_{_{j}}
(Y\!_{1},Y\!_{2};{\bf G}^{\prime})]
\longrightarrow}\hskip6cm
\nonumber\\
\hskip2.5cm\longrightarrow
\  {\hat{\mathbf{H}}}_{q+1}
        ((\!X\!_{_1},X\!_{_2})\
\hat{\ast} \ (Y\!_{_1},Y\!_{2}); {\bf
G}\otimes\!\!_{_{\mathbf{R}}}{\bf G}^{\prime})
\longrightarrow
\hskip4.5cm{{\hbox{\rm({\bf 3}\label{EqP14:3})}}}\!\!
\nonumber
\\
\hskip3cm\longrightarrow
\bigoplus_{i+j=q-1}
\ \ \hbox{\rm   Tor}_1^{\mathbf{R}}\bigl({\hat{\mathbf{H}}}_{_{i}}
           (X\!_{_1},X\!_{2};{\bf G}),
 {\hat{\mathbf{H}}}{_{j}}
      (Y\!\!_{_1},Y\!_{2};
{\bf G}^{\prime})\bigr)
\longrightarrow\ {\bf0}
\hskip1.5cm\Box
\nonumber
\end{array}$

\smallskip
Next theorem shows that also the local augmental homology for products and joins can be calculated through Theorem~\ref{Th04} once Proposition~\ref{Prop:000} p.~\pageref{Prop:000} (below) has been established.

\smallskip\noindent
{\bf Theorem 3.6.}\label{Th3.6} (See p.~\pageref{Th06}) {\it 
If $(t_{1},\widehat{x\ast y},t_{2}):=\{(x,y,t)\ {\mid}\ 0<t_{1}\le t\le t_{1}<1\}$
for polytopes $X$, $Y$ and $x\in X$, $y\in Y$ then,
{\rm (For definition of $\setminus_{o}$ see p.~\pageref{DefSetminus}. For missing proofs, see \cite{Fo1}.)}
$$
\mathbf{i. }\left\{\begin{array}{lcl}
\hat{\mathbf{H}}_{q+1}({X}\ast{Y},{X}\ast{Y}\setminus_{o}(x ,y ,t);\mathbf{G})\cong
\hat{\mathbf{H}}_{q+1}({X}\ast{Y},{X}\ast{Y}\setminus_{o}(t_{1},\widehat{x\ast y},t_{2});\mathbf{G})\cong \\
\cong\hat{\mathbf{H}}_{q}({X}\times{Y},{X}\times{Y}\setminus_{o}(x,y);\mathbf{G})\cong{\hat{\mathbf{H}}_{q}((X ,X \setminus_{o}x )\times(Y,Y\setminus_{o}y);\mathbf{G})}\cong \\
\cong\left[{{\text{\rm Motivation:}}\atop{\text{\rm Th.~\ref{Th02} line four.}}}\right]\cong\hat{\mathbf{H}}_{q+1}((X,X\setminus_{o}x)\ast(Y,Y\setminus_{o} y );\mathbf{G}).
\end{array}\right.
$$

$$\hskip-3.2cm
\mathbf{ii. }\ \hat{\mathbf{H}}_{q+1}(X \ast Y, X\ast Y\setminus_{o}(y ,{0});\mathbf{G})\cong\hat{\mathbf{H}}_{q+1} ((X ,\emptyset) \ast (Y,Y\setminus_{o}y);\mathbf{G})
$$
and equivalently for the $(x ,{1})$-points from the join-definition.\\
All isomorphisms are induced by chain equivalences. \hskip5cm$\Box$
}

\smallskip\noindent
{\bf Proposition~3.7.} \label{Prop:3.7.} {\rm(See p.~\pageref{Prop:000})} {\it Let $\mathbf{G}$ be a $($unital$)$ module over a commutative ring {\bf A} with unit.
With $\alpha$ in the interior of $\sigma$, i.e. $\alpha\in\{\beta\in|\Sigma|\mid[\mbox{v}\in\sigma] \Longleftrightarrow [\beta(\mbox{v})\neq0]\}$ and $\alpha=\alpha_{0}$ if and only if
$\sigma=\emptyset_o$,
the following {\bf A}-module isomorphisms are all induced by
chain $($homotopy$)$ equivalences. {\rm (For definition of $\setminus_{o},\ {\mathrm{Lk}}$ and $\mathrm{cost}$ see p.~\pageref{DefSetminus}. )}}
$$\hat{\mathbf{H}}_{i-\#\sigma}(\mathrm{ Lk}_{\Sigma}\sigma;\mathbf{G})\cong
\hat{\mathbf{H}}_{i}(\Sigma,\mathrm{ cost}_{\Sigma}\sigma;\mathbf{G})\cong\hat{\mathbf{ H}}_{i}(\vert\Sigma\vert,\vert\mathrm{ cost}_{\Sigma}\sigma\vert;
\mathbf{G})\cong\hat{\mathbf{ H}}_{i}(|\Sigma|,|\Sigma|\setminus{_o}\alpha;\mathbf{G}).\hskip0.5cm\Box$$

\newpage
\section{Stanley-Reisner rings  (Binary operations)}\label{SecP1:1}

\begin{definition}\label{simplDef00} An (augmented abstract) {\it simplicial complex} $\Sigma$ on a {\it vertex set} $V_\Sigma$ is a collection (possibly empty)
of finite subsets $\sigma$, the {\it simplices}, of $V_\Sigma$ satisfying:
$$\begin{array}{ll}
\mbox{({\bf a})}\text{ If } v\in V_{\Sigma}, \text{ then } \{v\}\in
\Sigma.\hskip3cm
\\
\mbox{({\bf b})} \text{ If } \sigma\in \Sigma \text{ and }
 \tau\subset \sigma, \text{ then } \tau\in \Sigma.
 \end{array}
$$
Define the {\it simplicial join} $\Sigma_1\ast\Sigma_2$, of two simplicial complexes $\Sigma_1$ and $\Sigma_2$ with \underbar{disjoint} vertex sets, to be:
$$
\Sigma_1\ast\Sigma_2:=\{\sigma_1\cup\sigma_2 \mid
\sigma_i\in\Sigma_i \ (i=1,2) \}.
$$
\end{definition}
Three levels can be identified: a {\it simplicial complex} is a set of{ \it simplices}, which are {\it finite} sets of {\it vertices}.
Denote the number of vertices in a simplex $\sigma\mbox{ by } \#\sigma = card(\sigma)$.

Put $ q := \dim\sigma := \#\sigma-1$. Then $\sigma$ is said to be a q-{\it face} or a q-{\it simplex} of $\Sigma$ and $$\dim\Sigma:=sup\{\dim(\sigma) | \sigma\in\Sigma\}.$$
Writing $\emptyset_o$ when using the empty set $\emptyset$ as a simplex, $\dim(\{\emptyset_o\})=-1$, while
$\dim(\emptyset)=-\infty$.

So, the empty set $\emptyset$ plays a dual role; it is both a $(-1)$-dimensional simplex and a $(-\infty)$-dimensional simplicial complex. The dimension formula for the join reads,
$$\dim(\Sigma_1\ast\Sigma_2)=\dim\Sigma_1+\dim\Sigma_2 +1,$$ implying that a join unit must be of dimension $(-1)$, as is the case for $\{\emptyset_o\}$.

$\emptyset$ is the unique zero-object with respect to join, i.e.,
$$\Sigma\ast\emptyset=\emptyset\ast\Sigma=\emptyset,
\mbox{ while } \{\emptyset_o\} \mbox{ is the join-unit,}$$
$$\Sigma\ast\{\emptyset_o\}=\{\emptyset_o\}\ast\Sigma=\Sigma.$$

\medskip\noindent
\begin{definition}\label{barDef} Let a bar over a set-symbol denote {\it the set of all \underbar{finite} subsets} of that set. This operation always results in a simplicial complex - known as the {\it the full simplicial complex} on that particular set.
Let the {\it boundary} of a set be ``the set of all its finite \underbar{proper} subsets" - indicated by a dot over the set-symbol.
\end{definition}
So, for any simplex:
$$
\bar\sigma:=\{\tau \mid \tau{\subset}\sigma\} \mbox{ and } {\dot{\sigma}}:=\bar{\sigma}{\setminus
\{\sigma\}}.$$
$$\text{ E.g., }\bar{\{v\}}=\{\emptyset_o,{\{v\}}\},\ {\bar{\emptyset}}_o =\{\emptyset_o\}=\dot{\{v\}}\text{ and }{\dot{\emptyset}}_o=\emptyset.
$$
\begin{definition}\label{defNonSimp} A finite subset
$s\subset W \supset V_\Delta$
is said to be a {\it non-simplex} (with respect to the ``universe" $W$) of a
simplicial complex $\Delta$ denoted $s\propto\Delta$ if: $s\not\in\Delta \mbox{ but }\dot{s}\subset\Delta.$
\end{definition}

\begin{definition} For a set $\delta=\{v_{i_1}\dots v_{i_k}\} \mbox{ let } m_{\delta} $  be the square-free monic monomial
$$m_{\delta}:=1_{\mathbf{A}}\cdot v_{i_1}\cdot v_{i_2}\cdot\dots\cdot v_{i_k}\in \mathbf{A}[W]$$
where $\mathbf{A}[W]$ is the graded polynomial algebra on the variable set $W$ over the commutative ring $\mathbf{A}$ with unit $1_{\mathbf{A}}$.
So, in particular, $m_{\emptyset_o}=1_{\bf A}.$
\end{definition}
The generating set for the ideal ${\mathbf{I}}_\Delta$ below is minimal, while in the
\htmladdnormallink{\it traditional definition }
{http://en.wikipedia.org/wiki/Stanley-Reisner_ring}
 of Stanley-Reisner rings, ${\mathbf{I}}_\Delta$ is generated by all square-free monomials
$v_{i_1}\cdot v_{i_2}\cdot\dots\cdot v_{i_k}\in \mathbf{A}[W]$ such that
$\{v_{i_1}, v_{i_2},\dots, v_{i_k}\}\notin\Delta,$ which is formally simpler but less rooted in logics and less explicit about the impact of different binary operations on simplicial complexes.
\begin{definition}\label{defStR-ring} Let
$$\mathbf{A} \langle \Delta \rangle := \mathbf{A}[W]/{\mathbf{I}}_{\Delta},$$
where  ${{\mathbf{I}}}_{\Delta}$ is the ideal generated by the monomials induced by the non-simplices of $\Delta,$ i.e. by $\{m_\delta  \mid \delta\propto\Delta \}.$
$\mathbf{A} \langle \Delta \rangle$ is called the {\it face ring} or the
{\it Stanley-Reisner $($St-Re$)$ ring} of $\Delta$ over $\mathbf{A}$.
Frequently $\mathbf{A}={\mathbb{Z}},$ the integers, or $\mathbf{A}={\mathbf{k}} ,$ a field, e.g. a prime field, the reals or the rational numbers.
\end{definition}

\noindent
Example {\bf ii} below shows that:
$ \mathbf{A} \langle\{\emptyset_o\}\rangle \cong \mathbf{A} = \text{The base-ring}$
$\ne{\mathbf{0}} =\mbox{The trivial ring}=\nobreak \mathbf{A}\nobreak\langle\emptyset\rangle.$

\medskip\noindent
{\bf Examples:}

\smallskip\noindent
{\bf i.}\label{Ex1:01} The choice of the universe $\mathbf{W}$ isn't all that critical, since
$$
\mathbf{A} \langle \Delta \rangle
\cong {{\mathbf{A}[V_\Delta]}\over(\{m_\delta\in\mathbf{ A}[V_\Delta]
\mid \delta\propto\Delta\})}, \text{ if } \Delta\not= \emptyset, \{{\emptyset}_o\}.
$$
{\bf ii.}\label{Ex1:02} If  $\Delta=\emptyset,$ then the set of non-simplices equals $\{\emptyset_o\}$, since $\emptyset_o\not\in\emptyset$ and

$\dot{\emptyset}_o={\bar{\emptyset}_o}^{((\dim\emptyset_o)-1)}=
{\{\emptyset_o\}}^{(-2)}=\emptyset\subset\emptyset\Longrightarrow$
$
{\mathbf{A}}\langle\emptyset\rangle=\mathbf{0}=\text{the trivial ring, since }
m_{\emptyset_{o}}=1_{\mathbf{A}}.$

Since $\emptyset \in \Delta$ if $\Delta\neq\emptyset,$
$\{v\}$ is a non-simplex of
$\Delta$ for every
$v\in\mathbf{W}\setminus \mathit{V}_{\Delta}$ i.e.,
$$[v\not\in\mathit{ V}_{\Delta}\neq\emptyset
{{]}{\Longleftrightarrow}[}\{v\}\propto\Delta\neq\emptyset].$$
\indent So
$\mathbf{A} \langle \{\emptyset_{o}\} \rangle= \mathbf{ A}$ since now $\{\delta\mid{\delta}\propto\Delta\}= \mathbf{W}.$

\medskip\noindent
{\bf iii.}\label{Ex1:03}
$
\mathbf{k}\ \langle
\Delta_1 \ast\Delta_2
\rangle
\cong \mathbf{k} \langle
\Delta_{1} {\rangle}
{\otimes}_\mathbf{k} \mathbf{k} {\langle} \Delta_{2} {\rangle}
$ (R. Fröberg, 1988.) and
$$
\mathbf{I}_{\Delta_{1}\ast\Delta_{2}}
=
( \{m_{\delta}\mid[
{\delta}\propto\Delta_1 \lor {\delta}\propto\Delta_2] \land [\delta\notin\Delta_1 \ast\Delta_2]\}).
$$

$(\cdot):=$the ideal generated by $\cdot$, $\land:=$``and" while  $\lor:=$``or".

\medskip\noindent
{\bf iv.} If $\Delta_i\not=\emptyset \text{ for } i=1,2,$ it is well known that,
$$\begin{array}{ll}
\mbox{\bf a. }
\mathbf{I}_{{\Delta}_{1}\cup {\Delta}_{2}}
=
\mathbf{I}{_{{\Delta}_{1}}}\cap  \mathbf{I}_{\Delta_{2}}
=
(\{m=\mathrm{Lcm}(m_{\delta_1},m_{\delta_2})\mid
\delta_i{\propto\Delta_i} \text{ for } i=1,2 \});\medskip
\\
 \mbox{\bf b. }
\mathbf{I}_{{\Delta}_{1}\cap{\Delta}_{2}}
=
\mathbf{I}_{{\Delta}_{1}}+ \mathbf{I}_{{\Delta}_{2}}
= (\{m_{{\delta}}\mid\delta\propto\Delta_1\lor \delta\propto\Delta_2\} ) \text{ in } \mathbf{A}[\mathbf{W}].
\end{array}
$$
$
\mathbf{I}_{\Delta_1}\cap \mathbf{I}_{\Delta_2}
\text{ and } \mathbf{I}_{\Delta_1}+ \mathbf{I}_{\Delta_2}
$
are generated by a set (no restrictions on its cardinality) of
square-free monomials if both $\mathbf{I}_{\Delta_1} \text{ and } \mathbf{ I}_{\Delta_2}$ are square-free.

\begin{definition}\label{defOrderSimpC} An {\it ordered simplicial complex}  $\Delta$ is a simplicial complex equipped with a partial order on $\mathbf{V}_{\Delta}$ that induces a linear order on the simplices.
\end{definition}

\begin{definition}\label{defOrderSimpProd} Given two ordered simplicial complexes $\Delta^\prime$, $\Delta^{\prime\prime}$
with vertex sets
$
V_{\Delta^\prime}:=
\{v_{1}^\prime,\dots, v_{a}^\prime \}
$
and
$
V_{\Delta^{\prime\prime}}:=
\{v_{1}^{\prime\prime},\dots, v_{b}^{\prime\prime}\}
$
resp., where all the
vertices belong to a common ``universe" $\mathbf{W}$.
\cite{E&S} defines in Def. 8.8 p. 67
{\it the ordered simplicial Cartesian Product}, $\Delta^\prime\times\Delta^{\prime\prime},$
of $\Delta^\prime$ and $\Delta^{\prime\prime},$
and shows in Lemma 8.9 p. 68 that
$
\Delta^\prime\times\Delta^{\prime\prime}
$
triangulates
$
|\Delta^\prime|\times|\Delta^{\prime\prime}|
$ - the topological cartesian product of the realizations of the involved complexes. Now, put,
$$
V_{\Delta^\prime\times\Delta^{\prime\prime}}:=
\{(v_{1}^\prime,v_{1}^{\prime\prime}),\dots,
(v_{a}^\prime,v_{b}^{\prime\prime})\} \text{ - the vertex set of } \Delta^\prime\times\Delta^{\prime\prime}
\text{ and } w_{i,j}:=(v_{i}^\prime,v_{j}^{\prime\prime}).
$$
Simplices in $\Delta^\prime\times\Delta^{\prime\prime}$ are
sets
$
\{w_{i_0,j_0}, w_{i_1,j_1},\dots, w_{i_k,j_k}\}, \text{ with } w_{i_{s},j_{s}}\not= w_{{i_{s+1}},{j_{s+1}}} \text{ and}
$
$$
v_{i_0}^\prime\le v_{i_1}^\prime\le\dots \le v_{i_k}^\prime \ \
(v_{j_0}^{\prime\prime}\le v_{j_1}^{\prime\prime}\le\ldots
\le v_{j_k}^{\prime\prime})
\text{  where } v_{i_0}^\prime, v_{i_1}^\prime,\dots, v_{i_k}^\prime\ \
(v_{j_0}^{\prime\prime}, v_{j_1}^{\prime\prime},\dots,v_{j_k}^{\prime\prime})$$
is a sequence of vertices, with repetitions possible, which constitutes a simplex in
$
\Delta^\prime\ \  (\Delta^{\prime\prime}).
$
\end{definition}
In particular,
$
\Delta^\prime\times\Delta^{\prime\prime}
$
is an ordered simplicial complex with respect to the product order.

$
\dim(\Delta^\prime\times\Delta^{\prime\prime}) = \dim\Delta^\prime + \dim\Delta^{\prime\prime}
$
except for the following products:
$$
\Delta \times \emptyset= \emptyset \times \Delta= \emptyset,
\text{ while } \Delta \times \{\emptyset_o\}= \{\emptyset_o\} \times \Delta = \{\emptyset_o\} \text{ if } \Delta\not=\emptyset.
$$

\medskip\noindent
{\bf Example v.} (\cite{Fo2} p.~72)
$C^\prime\cup D$
below is a reduced (Gröbner) basis, for
$\mathbf{I}$ in
$$
{k} {\langle}
\Delta_{1}\times\Delta_{2}
{\rangle}
\cong{k}[\text{V}_{\Delta_1}\times
\text{V}_{\Delta_2}] / {\mathbf{I}},
$$
i.e.,
$
C^\prime \cup D=\{m_\delta\mid  \delta
\propto {\Delta_1
\times {\Delta_2}\}}
$
with $C^\prime$ and $D$ defined as follows,
$$C^\prime:=\{ w_{{\lambda,\mu}}w_{{\nu,\xi}}|
\lambda<\nu\land \mu>\xi\}, \text{ where }
w_{{\lambda,\mu}}:= (v_{{\lambda}},v_{{\mu}}), \text{ with }
v_{{\lambda}}\in
\text{V}_{{\Delta_{1}}}\text{ and }
v_{\mu}\in\text{V}_{{\Delta_2}}.$$

The subindices reflect the assumed linear ordering on the factor simplices.\\
Now, letting $p_i$ with a bar over it denote the projection down onto the i:th factor;
$$
D:=\Bigl\{{w} =
w_{{\lambda_{1}},{\mu_{1}}}
\cdot\dots\cdot
w_{{{\lambda_{{k}}},{\mu_{_{^{k}}}}}} \mid\Bigl[
\bigl[\bigl[\{\overline{p_1}({w})\} \propto {\Delta_1}\bigr]
\land\bigl[\{\overline{p_2}({w})\}
\in{\Delta_2} \bigr] \land
\Bigl[{{\lambda_1<\dots<\lambda_k}\atop{\mu_1\le\dots\le \mu_k}}\Bigr]\bigr]\lor
$$
$$
\lor
\bigl[\bigl[ \{\overline{p_1}({w})\} \in{\Delta_1}
\bigr]
\land
\bigl[\{\overline{p_2}({w})\}\propto {\Delta_2}\bigr]
\land\Bigl[{{\lambda_1\le\dots\le \lambda_k} \atop {\mu_1<\dots<\mu_k}}\Bigr]\bigr]
$$
$$
\lor
\bigl[\bigl[\{\overline{p_1}({w})\}\propto{\Delta_1}\bigr]
\land  \bigl[\{\overline{p_2}({w})\}\propto{\Delta_2}\bigr] \land
\Bigl[{{\lambda_1<\dots<\lambda_k}\atop{\mu_1<\dots<\mu_k}}\Bigr]\bigr]\Bigr]\Bigr\}.
$$
Any Stanley-Reisner ring is a discrete Hodge algebra, as defined in \cite{Br} \S7.1.

The identification
$
v_{\lambda}\otimes v_{\mu}\leftrightarrow\
(v_{\lambda}, v_{\mu})
$ and the vertex orderings above,
gives the following graded {\bf k}-algebra
isomorphism of degree zero;
$$
{k} \langle
\Delta_{1}\times\Delta_{2}
\rangle
\cong
{k} \langle
\Delta_{1} \rangle
{\bar\otimes}
{k} \langle
\Delta_{2}\rangle,
$$
where the r.h.s. is an example of a generator-order sensitive Segre product.

\newpage

\section{Augmental homology}
Augmental homology, denoted $\hat{\mathbf{H}}_{\ast}$ wether the simplicial homology functor or the singular homology functor is concerned, is a homology theory fulfilling all the Eilenberg-Steenrod axioms and unifies the classical relative homology functor and the ad-hoc invented functor known as the classical reduced homology functor.

There is also an augmental cohomology functor which relates to augmental homology in the same way as classical relative cohomology relates to classical relative homology.

The underlying technicalities are analogous to those in the classical introduction of simplicial homology and by just hanging on to the $\{\emptyset_o\}$-augmented chains, also when defining relative chains $_o$, one gets the ''relative simplicial augmental homology functor for $\mathcal{K}_o$-pairs'' (defined below), denoted $\hat{\mathbf{H}}_{\ast}$.

The same comment applies to singular homology, via the introduction of the $(-1)$-dimensional topological space $\{\wp\}$ explained below. An introduction is also laid out in Wikipedia ``Homology".

\subsection{Background}
Abstract simplicial complexes can be either ``classical" or ``augmented". Some 40 years ago combinatorialists begun to systematically use Stanley-Reisner rings, also known as ''face rings'', and therethrough commutative algebra to solve combinatorial problems. This process opened up with an algebraization of the classical definition of abstract simplicial complexes which resulted in the augmental version, which in turn required a redefinition of the
simplicial join, making it ``tensor-like".

E.g. the empty simplicial complex switched from being the unit object w.r.t. to simplicial join to instead become its zero object. The role as the join unit was taken over by the new vertex-free simplicial complex $\{\emptyset\}$, containing nothing but the empty simplex.

The following definition gives; \\(\mbox{\bf 1.}) ``the classical abstract simplicial complexes"
if the word ``non-empty" is included and \\
(\mbox{\bf 2.})
the contemporary definition, i.e., that of ``augmented abstract simplicial complexes" if non-empty is excluded.

\smallskip\noindent
{\bf Definition.} A ''simplicial complex'' is a set S of (\underbar{non-empty}) finite sets closed under the formation of subsets.
$$\mbox{Formally: }  \mbox{\bf 1. } (\sigma\in \Sigma)\land(\emptyset\neq\tau \subset \sigma)\Longrightarrow\tau\in\Sigma \mbox{ resp., \mbox{\bf 2.} }(\sigma\in\Sigma)\land(\tau\subset\sigma)\Longrightarrow\tau\in\Sigma.$$
The Stanley-Reisner ring of the empty simplicial complex is the trivial ring with the zero-element as its only member, while the Stanley-Reisner ring of $\{\emptyset\}$ is the base ring. The trivial ring, $\mathbf{0}$, is a subring of every ring. The (new) augmented simplex $\emptyset$ generates a subcomplex $\{\emptyset\}$ of every non-empty augmented abstract simplicial complex.

The category of topological spaces and continuous functions were left unattended, which for instance implied that there were no topological space that could serve as the realization (see below) of the new simplicial complex $\{\emptyset\}$. Suddenly there were two vertex free simplicial complexes, $\emptyset$ and $\{\emptyset\}$, but still only one point-free topological space, i.e. $\emptyset$.

The empty set $\emptyset$ plays a dual role in the contemporary category of simplicial complexes - it is both a $(-1)$-dimensional simplex and a $(-\infty)$-dimensional simplicial complex.
Let $\emptyset_{o}$ denote the empty set $\emptyset$ when regarded as a simplex and not as a simplicial complex.

\smallskip
\subsection{Definitions}
Since the presence of this empty simplex $\emptyset_{o}$ in every non-empty (augmented abstract) simplicial complex is the only thing that makes them different from the classical ones, it is straightforward to find useful functors between these two categories. The new simplicial complex $\{\emptyset_o\}$ will also serve as join-unit and as the (abstract simplicial) $(-1)$-standard simplex.

Let ${{\mathcal K}}$ be the classical category of simplicial
complexes and simplicial maps and let
${{\mathcal K}_{o}}$ be the contemporary category.
Define $\mathcal{E}_{o}:{\mathcal K}\rightarrow{{\mathcal
K}_{o}}$ to be the functor augmenting $\emptyset_{o}$,
as a simplex, to each classical simplicial complex.
$\mathcal{E}_{o}$ has an inverse $\mathcal{E}:{\mathcal
K}_{o}\rightarrow{\mathcal K}$ deleting $\emptyset_{o}.$

The dimension of any non-empty simplicial complex is a well-defined integer $\geqq-1$, which for the join of two simplicial complexes with disjoint vertex sets results in the following dimension formula:
$$\dim (\Sigma_1\ast\Sigma_2)=\dim\Sigma_1 + \dim\Sigma_2 + 1,$$
implying that the join unit must be of dimension $-1$, which by definition is the dimension of $\{\emptyset_{o}\}$.

The least one expect from a realization functor is that it preserves the dimension.
So, this new $(-1)$-dimensional simplicial complex $\{\emptyset_o\}$ then actually implies that there can be no faithful formal realization functor from the augmented  simplicial complexes down to the classical topological spaces, in contrast to the purely classical situation.
The classical category of topological spaces only allows the \underbar{two} augmented vertex-free simplicial complexes $\emptyset$ (($-\infty$)-dimensional) and $\{\emptyset_o\}$ (($-1$)-dimensional) to be collapsed into the \underbar{one and only} point-free topological space $\emptyset$, which prevents the algebraic structure inflated in the 1970th into the  augmental category of simplicial complexes to be faithfully carried over to the category of topological spaces.
The same anomalous result occurs when the realization of augmented simplicial complexes are filtered through the category of simplicial sets, as worked out in \cite{E&P}.

Augmenting this new ($-1$)-dimensional topological space $\{\wp\}$ into each classical topological space constitutes the Augmental category of topological spaces, which certainly enable a full-fledged realization functor that is faithful to the new definition of simplicial complexes and their joins.
Of course, also the classical homology theory for pairs of simplicial complexes needs to be modified accordingly to make it adjust to this new algebra-adjusted simplicial environment.
On the chain level, every non-empty simplicial chain-complex remains the same as those in the classical theory, except for a single generator in degree $-1$, accounting for the everywhere present new ($-1$)-dimensional \underbar{simplex} $\emptyset_o$.

Using classical homology and the above terminology, the augmental homology theory for ``simplicial pairs" from the new simplicial category is completely determined as follows:
$$\hat{\mathbf{H}}_{i}(\Sigma_{o1},\Sigma_{o2};\mathbf{G})=
\left\{\begin{array}{llll}
\mbox{H}_i(\mathcal{E}(\Sigma_{o1}),\mathcal{E}(\Sigma_{o2});\mathbf{G}) & \mbox{ if } \Sigma_{o2}\not=\emptyset\\
\tilde{\mbox{H}}_i(\mathcal{E}(\Sigma_{o1}),\mathbf{G}) & \mbox{ if }
\Sigma_{o1}\neq\{\emptyset_o\}, \emptyset \mbox{ and } \Sigma_{o2}= \emptyset\\
\left\{\begin{array}{ll}
\cong \mathbf{G}  & \mbox{ if } i=-1\\
=0 & \mbox{ if } i\neq -1
\end{array}\right.
& \mbox{ when } \Sigma_{o1}=\{\emptyset_o\} \mbox{ and }
\Sigma_{o2}= \emptyset\\
0 & \mbox{ for all } i \mbox{ when } \Sigma_{o1}=\Sigma_{o2}=
\emptyset.
\end{array}\right.
$$
${\mbox{H}}_\ast$  denotes the
classical simplicial homology functor with its reduced tilde-equipped companion.
Note that now there is no need for any ``reduced homology functor" to the left.

To perform an analogous algebraization of the category of topological spaces - add, using topological sum, to each classical topological space $X\in \mathcal{D}$ an external element $\wp$, resulting\nobreak\ in;
$$X_{\wp}:=X + \{\wp\}\in\mathcal{D}_{\wp}.$$
Finally, add the universal initial object $\emptyset$ and let, as for $\mathcal{E}_{o}$ above,
$${\mathcal{F}_{\wp}}:{\mathcal{D}}\rightarrow{{\mathcal{D}}_{\wp}};X \longrightarrow {X}_{\wp}:=X + \{\wp\}$$
be the functor adjoining (i.e. ``augmenting") to each classical topological space a new non-final element denoted $\wp$.
$\mathcal{F}_{\wp}$ has an inverse $\mathcal{F}:{\mathcal D}_{\wp}\longrightarrow{\mathcal D}$ deleting ${\wp}$.

This far everything seems to work smoothly, but to be able to construct a homology theory, the new ``standard simplices" and, as in this case, the new ''singular simplices" must be defined, but prior to that the simplicial maps and the continuous functions must be properly defined/identified. The latter identification is  performed simply by declaring that no essentially new functions are allowed, in the sense that e.g. the ``new" topological functions are all ``topological unions" of the form
$$f_\wp=f+\mbox{Id}_{\{\wp\}}:{X}_{\wp}:= X+\{\wp\}\longrightarrow {Y}_{\wp}:= Y+\{\wp\}.$$
Since $\emptyset$ is declared to be a universal initial object, any empty map $\emptyset\longrightarrow X_\wp$ must also be accepted as an ``augmental map".
The empty topological space $\emptyset$ has played a fix formalized role within the Eilenberg-Steenrod axiomatic formalism and it is compatible, through the realization functor, with its $(-\infty)$-dimensional counterpart $\emptyset$ in the contemporary category of simplicial complexes, where $\{\emptyset\}\ (\neq\emptyset)$ is a natural $(-1)$-dimensional object.

Single-space simplicial and singular homology theories are rather perifer - mostly pair-space theories are used and then always together with the ``convention";
$${\mbox{H}}_i(X;\mathbf{G}):={\mbox{H}}_i(X, \emptyset;\mathbf{G}),$$
which in fact is more than a mere convention in that it tacitly assumes the chain complex of $\emptyset$ to have nothing but the trivial group in each degree. This has its origin in \cite{E&S} p.~3 where the true convention ``{\sl... $(X,\emptyset)$ is usually abbreviated by $(X)$ or simply, $X$}" is found.

\smallskip\noindent
{\bf Lemma.}
(Analogously for augmental cohomology, by raising the subindex to an index.)
$$\hat{\mathbf{H}}_i(X_{\wp1},X_{\wp2};\mathbf{G})=
\left\{\begin{array}{llll}
{\mbox{H}}_{\mathit{i}}(\mathcal{F}(X_{\wp1}), \mathcal{F}(X_{\wp2});\mathbf{G}) & \mbox{ if } X_{\wp2}\neq \emptyset\\
\tilde{\mbox{H}}_{\mathit{i}}(\mathcal{F}(X_{\wp1});\mathbf{G})
& \mbox{ if } X_{\wp1}\neq \{\wp\}, \emptyset \mbox{ and }
X_{\wp2}=\emptyset\\
\left\{\begin{array}{ll}
\cong\mathbf{G} & \mbox{ if } {\mathit{i}}=-1 \\
=0 & \mbox{ if } {\mathit{i}}\neq-1
\end{array}\right.
& \mbox{ when } X_{\wp1}= \{\wp\} \mbox{ and } X_{\wp2}=\emptyset\\
0 &  \mbox{ for all } {\mathit{i}} \mbox{ when } X_{\wp1}=X_{\wp2}=\emptyset.
\end{array}\right.$$
Both the augmental simplicial, as well as the augmental singular pair-space homology functors are formal ''homology theories'' in the sense of Eilenberg-Steenrod \cite{E&S}. This, in particular, implies that relative
\htmladdnormallink{{\it Mayer-Vietoris sequences} }{http://en.wikipedia.org/wiki/Mayer-Vietoris_sequence} for excisive couples of pair spaces can be used without further motivation, as this is a direct formal consequence of the
\htmladdnormallink{{\it Eilenberg-Steenrod axioms.} }{http://en.wikipedia.org/wiki/Eilenberg-Steenrod_axioms}

\smallskip
\subsection{Augmental homology for joins and products of augmented topological spaces}
The definition of {\it the product of simplicial complexes} is given in def.~\ref{defOrderSimpProd} p.~\pageref{defOrderSimpProd} at the end of the chapter on
{\it Stanley-Reisner rings $($Binary operations$)$},
that of simplicial sets is given in \cite{F&P} and the one for topological spaces is found in any elementary book on general topology. \\
The join of two simplicial complexes with disjoint vertex sets is a very important binary operation (cp. Ex.~iii p.~\pageref{Ex1:03}), defined as follows:
$$\Delta_1\ast\Delta_2:=\{\delta \mid \delta=\delta_1\cup\delta_2,\ \delta_i\in\Delta_i \},$$
while the join-definition for simplicial sets and topological spaces is found in \cite{F&G} resp. \cite{Fo1}.
The join of two topological spaces $X$ and $Y$ can be defined as the double mapping cylinder of $X\times Y$ with respect to the two projection maps $p_X: X\times Y\longrightarrow X$ and $p_Y: X\times Y\longrightarrow Y$, which might appear somewhat cryptic.
The detailed definitions of the two topological join-versions in common use, $\ast$ and $\hat{\ast}$, is found in \cite{Fo1} p.~17 and their actions are homotopic.

The pair-products and pair-joins are defined as follows:
$$(X_1,X_2)\times_\cup(Y_1,Y_2):=(X_1\times Y_1,X_1\times Y_2\cup X_2\times Y_1)$$
and
$$(X_1,X_2)\ast_\cup(Y_1,Y_2):=(X_1\ast Y_1,X_1\ast Y_2\cup X_2\ast Y_1).$$
The topologizing of topological pair-joins has to be done with some care, cf. \cite{Fo1} pp.~19-20.

These are the versions that are in common use, usually with $\diamond := \diamond_\cup$, though it is
$$(X_1,X_2)\times_\cap(Y_1,Y_2):=(X_1\times Y_1,X_1\times Y_2\cap X_2\times Y_1)=(X_1\times Y_1,X_2\times Y_2)$$
that is the pair-categorical product.
Equivalent definitions for $\ast_\cap$ and for simplicial pairs. \\

The {\it new} object $\{\wp\}$ gives the classical Künneth formula ($\equiv$4:th line in Th.~\ref{Th01} below) additional strength but much of the classical beauty is lost - a loss which is regained in the join version, i.e. in Theorem~\ref{Th04} p.~\pageref{Th04}
({\bf PID} = Principal Ideal Domain),

\begin{theorem}\label{Th01}
If $\{X_{1}\times Y_{2}, X_{2}\times Y_{1}\}$ is an ''excisive couple" {\rm(Def. \cite{Sp} p.~188)}, $\mbox{q}\geqq0$, ${\mathbf{R}}$ a ${\mathbf{PID}}$, and assuming $\mbox{Tor}_1^{\mathbf{R}}(\mathbf{G},\mathbf{G}^\prime)=0$ for ${\mathbf{R}}$-modules $\mathbf{G}$ and $\mathbf{G}^\prime$, then;
$$
\hat{\mathbf{H}}_{q}((X_{1},X_{2})\times (Y_{1}, Y_{2});\mathbf{G}\ {\otimes}_{\mathbf{R}}\mathbf{G}^\prime)\cong
$$
$$\left\{\begin{array}{llll}
[\hat{\mathbf{H}}_i(X_{1};\mathbf{G})\otimes_{\mathbf{R}}\hat{\mathbf{H}}_j(Y_{1};{\mathbf{G}}^{\prime})]_q\oplus
(\hat{\mathbf{H}}_{q}(X_{1};\mathbf{G})\otimes_{\mathbf{R}} {\mathbf{G}}^\prime)\oplus(\mathbf{G}\otimes_{\mathbf{R}}\hat{\mathbf{H}}_q Y_{1};{\mathbf{G}}^{\prime})
\oplus{\mbox{T}}_1 & \mbox{if}\ {\mbox{C}}_1\\
{[\hat{\mathbf{H}}_i(X_{1};\mathbf{G})\otimes_{\mathbf{R}}\hat{\mathbf{H}}_j(Y_{1},Y_{2};{\mathbf{G}}^\prime)]}_q \oplus({\mathbf{G}}\otimes_{\mathbf{R}}\hat{\mathbf{H}}_{q}(Y_{1},Y_{2};{\mathbf{G}}^\prime))\oplus\mbox{T}_2& \mbox{if}\ \mbox{C}_2\\
{[\hat{\mathbf{H}}_i(X_{1},X_{2};\mathbf{G})\otimes_{\mathbf{R}}\hat{\mathbf{H}}_j(Y_{1};{\mathbf{G}}^\prime)]}_q \oplus(\hat{\mathbf{H}}_{q}(X_{1},X_{2};\mathbf{G}) \otimes_{\mathbf{R}}{\mathbf{G}}^\prime)\oplus\mbox{T}_3 & \mbox{if}\ \mbox{C}_3,\\
{[\hat{\mathbf{H}}_i(X_{1},X_{2};\mathbf{G})\otimes_{\mathbf{R}}
\hat{\mathbf{H}}_j(Y_{1},Y{2};\mathbf{G}^\prime) ]}_q\oplus\mbox{T}_4 & \mbox{if}\ \mbox{C}_4,
\end{array}\right.({\bf 1})$$
where the torsion terms, i.e. the $\mbox{T}$-terms, split as those ahead of them, e.g.,
$$
\mbox{ T}_{1}=[\mbox{Tor}_{1}^{\mathbf{R}}\bigl(
\hat{\mathbf{H}}_i(X_{1};\mathbf{G}), \hat{\mathbf{H}}_j(Y_{1}; \mathbf{G}^\prime)\big)]_{q-1}
\oplus\mbox{Tor}_{1}^{\mathbf{R}}
\bigl(\hat{\mathbf{H}}_{q-1}(X_{1};\mathbf{G}), \mathbf{G}^\prime\bigr)\oplus\mbox{ Tor}_{1}^{\mathbf{R}}\bigl(\mathbf{G},\hat{\mathbf{H}}_{q-1}(Y_{1};\mathbf{G}^\prime)\bigr),
$$
and
$$ \mbox{ T}_{4}=[\mbox{Tor}_{1}^{\mathbf{R}}( \hat{\mathbf{H}}_i(X_{1},X_{2};\mathbf{G}),
\hat{\mathbf{H}}_j(Y_{1},Y_{2}; \mathbf{G}^\prime))]_{q-1}.
$$

$\mbox{ C}_i,\ (i=1-4)$, are ``conditions" and should be interpreted as follows, resp.,
$$\left\{\begin{array}{llll}
\mathrm{C}_1:=X_{1}\times Y_{1}\neq\emptyset,\{\wp\}\mbox{ and }X_{2}=\emptyset=Y_{2},\\
\mathrm{C}_2:=X_{1}\times Y_{1}\neq \emptyset,\{\wp\} \mbox{ and } X_{2}=\emptyset\neq Y_{2},\\
\mathrm{C}_3:=X_{1}\times Y_{1}\neq \emptyset,\{\wp\} \mbox{ and } X_{2}\neq\emptyset=Y_{2},\\
\mathrm{C}_4:=X_{1}\times Y_{1} = \emptyset,\{\wp\} \mbox{ or }  X_{2}\neq \emptyset\neq Y_{2}.\end{array}\right.
$$
$[\dots]_q \text{ above, should be interpreted as }
\bigoplus_{i+j=q \text{ and } i,j\geq0}\dots  .$ {\rm (... as in \cite{Sp} p.~235 Th. 10.)}.
\end{theorem}

\medskip
The next theorem is a consequence of the relative \htmladdnormallink{{\it Mayer-Vietoris sequence} }{http://en.wikipedia.org/wiki/Mayer-Vietoris_sequence} (M-Vs) and, in this case, the splitting of it. Moreover, note that excisivity is not an issue in this join/product-version of the M-Vs. The theorem shows the strong connection between the augmental singular homology groups for ``{\sl products of pair spaces}" and those for ``{\sl joins of pair spaces}". For $t\in[0,1]$, the intervalls $t\geqq0.5 \mbox{ and } t\leqq0.5$ identifies two halves of the pair join of two pair spaces, regarded as a collection of curves - one for each point pair. The union of these two halves is that very join, while the intersection is the pair product.

\medskip\noindent
\begin{theorem}\label{Th02}
For $(X_{1},X_{2})\neq(\{\wp\},\emptyset)\neq(Y_{1},Y_{2}) \text{ and if } \mathbf{A}$ is a commutative ring with unit and if $\mathbf{G}$ is an $\mathbf{A}$-module, then;

\smallskip
$
\hat{\mathbf{H}}_{q}((X_{1},X_{2})\times(Y_{1},Y_{2});\mathbf{G})
\cong
$
$$\cong\hat{\mathbf{H}}_{q+1}((X_{1},X_{2})\ast(Y_{1},Y_{2});\mathbf{G})\oplus
\hat{\mathbf{H}}_{_{q}}((X_{1},X_{2})\ast(Y_{1},Y_{2})^{^{_{t\geqq0.5}}}{+}
(X_{1},X_{2})\ast(Y_{1},Y_{2})^{^{t\leqq0.5}};\mathbf{G})=
$$
$$
\cong\left\{\begin{array}{llll}
\hat{\mathbf{H}}_{q+1}(X_{1}\ast Y_{1};
\mathbf{G}) \oplus \hat{\mathbf{H}}_{q}(X_{1};\mathbf{G}) \oplus
\hat{\mathbf{H}}_{q}(Y_{1};\mathbf{G}) & \mbox{if}\ \ \ \mbox{ C}_1
\\
\hat{\mathbf{H}}_{q+1}((X_{1},\emptyset)\ast(Y_{1},Y_{2});\mathbf{G})\oplus
\hat{\mathbf{H}}_{q}(Y_{1},Y_{2};\mathbf{G}) & \mbox{if}\ \ \ \mbox{ C}_2
\\
\hat{\mathbf{H}}_{q+1}((X_{1},X_{2})\ast(Y_{1},\emptyset);\mathbf{G})\oplus
\hat{\mathbf{H}}_{q}(X_{1},X_{2};\mathbf{G}) & \mbox{if}\ \ \ \mbox{ C}_3
\\
\hat{\mathbf{H}}_{q+1}((X_{1},X_{2})\ast(Y_{1},Y_{2});\mathbf{G}) & \mbox{if}\ \ \
\mbox{ C}_4
\ \ \
\end{array}\right.
\hskip3cm(\mathbf{2}),
$$
where the $+$-sign above indicates ``addition" at the underlying chain-level and
$\ \mbox{ C}_i,\ (i=1-4)$, are the same ``conditions" as in Th.~\ref{Th01} above.
\end{theorem}

\begin{theorem}\label{Th03}
{\rm[The relative Eilenberg-Zilber theorem for topological join.]}\\
 For an excisive couple $\{X\ast Y_{2},X_{2} \ast  Y\}$ from the category of ordered couples $((X ,X_{2}),(Y,Y_{2}))$ of topological pairs$_{\wp};$
$$\mathbf{s}({\Delta^{\wp}}(X,X_{2})\otimes{\Delta^{\wp}} (Y,Y_{2}))\ \text{\sl is naturally
\htmladdnormallink{{\it chain $($homotopy$)$ equivalent} }{http://en.wikipedia.org/wiki/Homotopy_category_of_chain_complexes}to } \Delta^{\wp}((X,X_{2})\ast(Y,Y_{2})).$$
{\rm($\mathbf{s}$ stands for  suspension i.e. the suspended
chain equals the original one except that the dimension $i$ in the original chain becomes $i+1$ in the suspended chain. Th.~\ref{Th03} is the augmental join version of the classical Th. 9 in \cite{Sp} p.~234 for products.)}
\end{theorem}

The couple
$\{X_{1}\ast Y_{2}, X_{2}\ast Y_{1}\}$ is excisive if and only if
$\{X_{1}\times Y_{2}, X_{2}\times Y_{1}\}$ is excisive.

\begin{theorem}\label{Th04}
{\rm[The \htmladdnormallink{{\it Künneth formula} }
{http://en.wikipedia.org/wiki/K\%C3\%BCnneth_theorem} for Topological Joins; cp. \cite{Sp} p.~235.]}\\
If $\{X_{1}\ast Y_{2}, X_{2}\ \ast Y_{1}\}$ is an excisive couple in $X_{1}\ast Y_{1},$ ${\mathbf{R}}$ a ${\mathbf{PID}}$, $\mathbf{G} \mbox{ and } \mathbf{G}^\prime\
\mathbf{R}$-modules and $\mathrm{Tor}_1^{\mathbf{R}}(\mathbf{G},\mathbf{G}^{\prime})=\mbox{\bf0}.$
Then the functorial sequences below are $($non-naturally$)$ split exact;
\begin{eqnarray}
{{\bf0}\longrightarrow \bigoplus_{i+j=q} \ \ [{\hat{\mathbf{H}}}_{_{i}}
           (X\!_{_1},X\!_{_2};{\bf G})
\otimes\!\!_{_{\mathbf{R}}}
        {\hat{\mathbf{H}}}_{_{j}}
(Y\!_{1},Y\!_{2};{\bf G}^{\prime})]
\longrightarrow}\hskip6cm
\nonumber
\\
\longrightarrow
\  {\hat{\mathbf{H}}}_{q+1}
        ((\!X\!_{_1},X\!_{_2})\
\hat{\ast} \ (Y\!_{_1},Y\!_{2}); {\bf
G}\otimes\!\!_{_{\mathbf{R}}}{\bf G}^{\prime})
\longrightarrow
\hskip4.5cm{{\hbox{\rm({\bf 3}\label{EqP14:3})}}}\!\!
\nonumber
\\
\hskip3cm\longrightarrow
\bigoplus_{i+j=q-1}
\ \ \hbox{\rm   Tor}_1^{\mathbf{R}}\bigl({\hat{\mathbf{H}}}_{_{i}}
           (X\!_{_1},X\!_{2};{\bf G}),
 {\hat{\mathbf{H}}}{_{j}}
      (Y\!\!_{_1},Y\!_{2};
{\bf G}^{\prime})\bigr)
\longrightarrow\ {\bf0}
\hskip1.5cm\Box
\nonumber
\end{eqnarray}

{\rm(\cite{Sp} p.~247 Th.~11 indicates how to formulate the cohomology-analog of Theorem~\ref{Th04}.)}
\end{theorem}
$(X_{1},X_{2})=(\{\wp\},\emptyset)$ in the last
Theorem immediately gives the following useful theorem. Note that the pair-space $(\{\wp\},\emptyset)$ does not exist classically. Moreover, general formulas containing joins of pair-spaces can not live in the classical environment, making classical algebraic topology a relatively week tool compared to the augmental version presented here. This is further underlined by the fact that the join-operation is the functorial counterpart, under the Stanley-Reisner ring construction, of the tensor product within commutative algebra, see Example iii p.~\pageref{Ex1:03} in chapter 2 on  Stanley-Reisner rings.

\begin{theorem}\label{Th05}
{\rm[$\mbox{The Universal Coefficient Theorem for (co)}\hat{\mathbf{H}}omology$]}
$$\left\{\begin{array}{ll}
\hat{\mathbf{H}}_{i}(Y_{1},Y_{2};\mathbf{G})\cong
\hat{\mathbf{H}}_{i}(Y_{1},Y_{2};{\mathbf{R}}\otimes_{\mathbf{R}}\mathbf{G})\cong\\
\cong (\hat{\mathbf{H}}_{i}(Y_1,Y_2;{\mathbf{R}})\otimes_{\mathbf{R}}\mathbf{G})\oplus
\mbox{Tor}^{{\mathbf{R}}}_1 \big(\hat{\mathbf{H}}_{i-1} (Y_{1},Y_{2};{\mathbf{R}}),\mathbf{G}\big)
\end{array}\right.
\ \ \ \text{ for any {\bf R}-{\bf PID}-module {\bf G}.}
$$
If all $\hat{\mathbf{H}}_{{\ast}}(Y_{1},Y_{2};{\mathbf{R}})$ are of finite type or $\mathbf{G}$ is finitely generated, then;
$$
\hat{\mathbf{H}}^{i}(Y_{1},Y_{2};\mathbf{G})\cong\hat{\mathbf{H}}^{i} (Y_{1},Y_{2};\mathbf{R}\otimes_{\mathbf{R}}\mathbf{G})
\cong
(\hat{\mathbf{H}}^{i}(Y_{1},Y_{2};\mathbf{R}) \otimes_{\mathbf{R}}\mathbf{G})\oplus \mbox{Tor}^{\mathbf{R}}_1 (\hat{\mathbf{H}}^{{{i+1}}}(Y_{1},Y_{2};{\mathbf{R}}), \mathbf{G}).
$$
\end{theorem}

Next theorem shows that also the local augmental homology for products and joins can be calculated through Theorem~\ref{Th04} once Proposition~\ref{Prop:000} p.~\pageref{Prop:000} has been established.

\begin{theorem}\label{Th06}
If $(t_{1},\widehat{x\ast y},t_{2}):=\{(x,y,t)\ {\mid}\ 0<t_{1}\le t\le t_{1}<1\}$
for polytopes $X$, $Y$ and $x\in X$ and $y\in Y$ then, {\rm (For definition of $\setminus_{o}$ see p.~\pageref{DefSetminus}. )}
$$
\mathbf{i. }\left\{\begin{array}{lll}
\hat{\mathbf{H}}_{q+1}({X}\ast{Y},{X}\ast{Y}\setminus_{o}(x ,y ,t);\mathbf{G})\cong
\hat{\mathbf{H}}_{q+1}({X}\ast{Y},{X}\ast{Y}\setminus_{o}(t_{1},\widehat{x\ast y},t_{2});\mathbf{G})\cong \\
\cong\hat{\mathbf{H}}_{q}({X}\times{Y},{X}\times{Y}\setminus_{o}(x,y);\mathbf{G})\cong{\hat{\mathbf{H}}_{q}((X ,X \setminus_{o}x )\times(Y,Y\setminus_{o}y);\mathbf{G})}\cong \\
\cong\left[{{\text{\rm Motivation:}}\atop{\text{\rm Th.~\ref{Th02} line four.}}}\right]\cong\hat{\mathbf{H}}_{q+1}((X,X\setminus_{o}x)\ast(Y,Y\setminus_{o} y );\mathbf{G}).
\end{array}\right.
$$
$$\hskip-3.2cm
\mathbf{ii. }\ \hat{\mathbf{H}}_{q+1}(X \ast Y, X\ast Y\setminus_{o}(y ,{0});\mathbf{G})\cong\hat{\mathbf{H}}_{q+1} ((X ,\emptyset) \ast (Y,Y\setminus_{o}y);\mathbf{G})
$$
and equivalently for the $(x ,{1})$-points from the join-definition.\\
All isomorphisms are induced by chain equivalences. \hskip5cm$\Box$
\end{theorem} 

The proofs of the above theorems are all collected in \cite{Fo1}. \\
The Proposition below and the link-formula:
$$\mbox{Lk}_{\Sigma_1\ast\Sigma_2}(\sigma_1\cup\sigma_2)= (\mbox{Lk}_{\Sigma_1}\sigma_1)\ast (\mbox{Lk}_{\Sigma_2}\sigma_2),$$
offers a good platform for a straightforward proof of Theorem~\ref{Th06} in this special case of polytopes - and therefore also for \htmladdnormallink{{\it CW-complexes} }{http://en.wikipedia.org/wiki/CW_complex} since those have the homotopy type of the realization of a simplicial complex and since joins preserves homotopies and the singular homology functor isn't pair-homotopy sensitive.

\smallskip
\subsection{The ``point" of augmental structures}
There is a notable and maybe rather annoying effect of the algebraization, i.e. the shift from classical complexes/spaces to augmented complexes/spaces resp. concerning the normal perception of a ``point" in every-day-mathematics. In algebraic topology the concept of a ``point" is \underbar{formalized} and of fundamental importance, mainly due to its presence in the \htmladdnormallink{{\it dimension axiom} }{http://en.wikipedia.org/wiki/Eilenberg-Steenrod_axioms} among the Eilenberg-Steenrod axioms, which, on a formal level, is non-problematic to implement in the new categories of augmented simplicial complexes resp. augmented topological spaces.
In classical simplicial homology theory a point is a simplicial complex containing a single simplex that contains a single vertex, formally $\{\{{v}\}\}$. In the literature this \underbar{classical simplicial complex} $\{\{{v}\}\}$ is often denoted $\{{v}\}$ which formally denotes a \underbar{classical simplex}.

``Points" in a complex ${\Sigma}$ in the category of augmented simplicial complexes $\mathcal{K}_{o}$ now becomes $\{\emptyset,\{{v}\}\}$, which certainly does not match our intuitive perception of a point, but within the category of augmented simplicial complexes and under the Eilenberg-Steenrod's axiomatization of (co)homology theory,
$\{\emptyset,\{{v}\}\}$
is formally a typical ``point".

The Stanley-Reisner ring
(St-R ring, Def. \ref{defStR-ring} p.~\pageref{defStR-ring}) of $\emptyset$ is the trivial ring $\mathbf{0}$. In the light of the algebraization --- as awkward as it would be to define and work with rings with no trivial subring $\mathbf{0}$, as natural will it probably become to work with ``points" of the type $\{\emptyset,\{{v}\}\}$ in augmented complexes $\Sigma\in\mathcal{K}_{o}$, i.e. in the category of augmented simplicial complexes, while of the generic type $\{\wp,x\}$ in augmented spaces
$X_{\wp}\in\mathcal{D}_{\wp},$
the category of augmented topological spaces (= augmental category of topological spaces).

The notation ``$x\in X$" is just a short for the set-theoretical ``$\{x\}\subset X$", which in\nobreak\ turn, via $\mathcal{F}_{\wp}$, impose the notation ``$\{\wp,x\}\subset X_{\wp}$" for a {\sl point} in $X_{\wp}\in\mathcal{D}_{\wp} =$ the category of augmented topological spaces and continuous functions. Well, no harm comes out of adopting the convention that ``$\{\wp,x\}\subset X_{\wp}$" also in the augmental setting is denoted\nobreak\ $x\in X.$

\subsection{More on notations and realizations}
The next Proposition is an augmental version of results from 1984 due to H-G. Gräbe and J.R. Munkres. (Changing the subindex $i$ to an index gives the formula for augmental cohomology.)
It shows a very close homology connection between the local structure of a simplicial complex and that of its polytope. This close connection is in itself a strong motivation for the introduction of the topological (-1)-dimensional object $\{\wp\}$, which, due to its relation to $\{\emptyset_o\}$, imposes the following definition of an extended ``{\sl setminus}" ``${\setminus_o}$" in ${{\mathcal D}_{\wp}}$, where the identification $x\leftrightarrow\{x,\wp\}$, modeled after the classical $x\leftrightarrow\{x\}$, is used. The common practise to write $x\in X$ for the classical point-set-theoretical $\{x\}\subset X$ motivates the above notation $x\in X$ also for $\{\wp,x\}\subset X_{\wp}.$

``$\setminus$" denotes the {\it classical setminus} also known in set theory as the relative complement.

To avoid to ``accidentally" dropp out of the new category, the definition of an extended ``setminus" ($\setminus_o$) needs some extra attention:

\smallskip\noindent
{\bf Definition.}\label{DefSetminus}
$$X_{\wp1}\setminus_o X_{\wp2}:=
\left\{\begin{array}{ll}
\emptyset & \mbox{ if }X_{\wp1}=\emptyset, X_{\wp1}\subsetneqq
X_{\wp2}
\mbox{ or }X_{\wp2}=\{{\wp}\}\\
\mathcal{F}_{\wp}({\mathcal{F}}(X_{\wp1})
{\ }\setminus{\ }
{\mathcal{F}}(X_{\wp2})) & \mbox{ else}.
\end{array}\right.
$$

In particular, $X_{\wp1}\setminus_o{x}=\emptyset\ \ \mbox{if}\ \ x=\wp$.

The {\sl link}, $\mbox{Lk}_{\Sigma}\sigma$, of a simplex $\sigma$ with respect to a simplicial complex $\Sigma$ is defined through: \label{DefLk}
$$
{\mathrm{Lk}}_{\Sigma}\sigma:= \{\tau\in \Sigma\mid
[\sigma\cap \tau =\emptyset]\land [\sigma\cup \tau \in
\Sigma]\} $$
Note that $\mbox{Lk}_{_{\Sigma}}\emptyset_{_{^o}}=\Sigma$ and this very link is sometimes called ``the Missing Link", since it does not exist in the classical setting.
$$\dim\mbox{Lk}_{\Sigma}\sigma=\dim\Sigma-\#\sigma,$$
where $\#$ stands for cardinality, i.e., $\#\sigma$ denotes ``the number of vertices in the simplex $\sigma$".

The {\sl contrastar} of $\sigma$ w.r.t. $\Sigma$ is defined as:
$$\mathrm{cost}_{\Sigma}\sigma:= \{\tau\in\Sigma|
\tau\nsupseteqq\sigma\}.$$
$\mbox{  So, cost}_{{\Sigma}}\emptyset_{o}=\emptyset
\mbox{ and }
\mbox{cost}_{{\Sigma}}\sigma=\Sigma \mbox{ if and only if } \sigma\not\in\Sigma.$
(In classical literature, the ``$\mbox{cost}_{{\Sigma}}\sigma$" is known as the ``complement of $\Sigma$ w.r.t. $\sigma$", e.g., see \cite{E&S} p.~74 exercise E.)

Now, with the $(-1)$-dimensional topological space $\{\wp\}$ in place, the new definition of the realization $\vert\Sigma\vert$ of an arbitrary simplicial complex $\Sigma$ is carried out by just extending the range space (a function space) of Spanier's classical covariant realisation-functor from \cite{Sp} p.~110, with an additional ``coordinate function" $\alpha_0$ where
$\alpha_0(v)\equiv 0 \mbox{ for all}\mbox{ v}\in \mbox{V}_{{\Sigma}}, \mbox{ where V}$ is the vertex set.
In the resulting function space, $\{\alpha_0\}$ now serves as the $(-1)$-dimensional topological join unit $\{\wp\}$, i.e.,  $\vert\{\emptyset_{o}\}\vert=\{\alpha_0\}$.
For details, see \cite{Fo1} p.~13.

Since the barycentric realization of the ($-1$)-dimensional simplicial complex $\{\emptyset_o\}$ is the ($-1$)-dimensional topological space $\{\alpha_0\}$ where $\alpha_0$ is the unique coordinate function that maps any vertex to $0$, the distance from $\alpha_0$ to any other point is constant, i.e., it is exactly $1$ in the barycentric metric --- a truly outstanding property.
\begin{proposition} \label{Prop:000} {\rm(For proof, see \cite{Fo1} p.~24.)} Let $\mathbf{G}$ be a $($unital$)$ module over a commutative ring {\bf A} with unit.
With $\alpha$ in the interior of $\sigma$, i.e. $\alpha\in\{\beta\in|\Sigma|\mid[\mbox{v}\in\sigma] \Longleftrightarrow [\beta(\mbox{v})\neq0]\}$ and $\alpha=\alpha_{0}$ if and only if
$\sigma=\emptyset_o$,
the following {\bf A}-module isomorphisms are all induced by
chain $($homotopy$)$ equivalences.
$$\hat{\mathbf{H}}_{i-\#\sigma}(\mathrm{ Lk}_{\Sigma}\sigma;\mathbf{G})\cong
\hat{\mathbf{H}}_{i}(\Sigma,\mathrm{ cost}_{\Sigma}\sigma;\mathbf{G})\cong\hat{\mathbf{ H}}_{i}(\vert\Sigma\vert,\vert\mathrm{ cost}_{\Sigma}\sigma\vert;
\mathbf{G})\cong\hat{\mathbf{ H}}_{i}(|\Sigma|,|\Sigma|\setminus{_o}\alpha;\mathbf{G}).$$
\end{proposition}
Since the realization of a simplicial complex is a compactly generated normal space it is convenient to assume that the  realization-functor targets the category of compactly generated \underbar{Hausdorff} spaces, where the concept ``compactly generated" is unambiguous, compare \cite{May} Ch.~5 p.~39ff. This topology-enlargement to the compactly generated topology is rarely needed in practice and is in any case harmless as it, e.g., does not effekt the calculations of homology groups. The convenience depends on the fact that both the topological join- and the product-topology of the realisations of two simplicial complexes otherwise need not in general be compactly generated. With this in mind, the formulas for the realization of a product $(\times)$ or a join $(\ast)$ of simplicial complexes becomes very simple since the simplicial joins and products just turn into topological joins and products resp.

The Milnor realization $\hat{\vert}\Xi\hat{\vert}$
of any simplicial set $\Xi$
is triangulable by \cite{F&P} p.~209 Cor.\ 4.6.12.

E.g., the augmental singular complex
$\Delta^{\wp}(X)$ with respect to any topological space $X$, is a simplicial set and, cf. \cite{Mi} p.~362 Theorem 4, the map
$
j:\hat{\vert}\Delta^{\wp}(X)\hat{\vert}\rightarrow X
$
is a weak  \htmladdnormallink{{\it homotopy equivalence} }{http://en.wikipedia.org/wiki/Homotopy_category_of_chain_complexes} i.e. induces homotopy isomorphisms, and $j$ is a true
homotopy equivalence if $X$ is of CW-type, e.g. if $X$ is a polytope, cf. \cite{F&P} pp.~76-77, 170, 189ff, 221-2.

Any (augmental) simplicial complex $\Sigma$ can be regarded as an ordered simplicial complex.
Any ordered simplicial complex $\Sigma$ also possess the structure of a simplicial set. If the latter structure is at focus the symbol for the simplicial complex is equipped with a hat, as is any symbol that relates to simplicial sets in the following formulas. Product respectively join signs within the category of compactly generated spaces are below equipped with a bar over it. Now, \cite{F&P} and \cite{F&G}, resp., gives the following homeomorphism for complexes $\Sigma_1, \Sigma_2$:
$$
\vert\Sigma_{1}\vert{\bar\times}\vert\Sigma_{2}\vert\simeq\hat{\vert}
\hat\Sigma_{1}\hat{\vert}\bar\times\hat{\vert}\hat\Sigma_{2}\hat{\vert}
\simeq
\hat{\vert}\hat\Sigma_{1}\hat\times\hat\Sigma_{2}\hat{\vert}
\simeq
\vert\Sigma_{1}\times\Sigma_{2}\vert$$
and
$$\vert\Sigma_{1}\vert{\bar\ast}\vert\Sigma_{2}\vert
\simeq
\hat{\vert}\hat\Sigma_{1}\hat{\vert}\bar\ast
\hat{\vert}\hat\Sigma_{2}\hat{\vert}
\simeq
\hat{\vert}\hat\Sigma_{1}\hat\ast\hat\Sigma_{2}\hat{\vert}
\simeq
\vert\Sigma_{1}\ast\Sigma_{2}\vert.
$$
The realization functor and the functor induced by
$
j:\hat{\vert}\Delta^{\wp}(X)\hat{\vert}\rightarrow X$ are adjoint functors in the sense of Kan \cite{KanAdjFunc}, making the structures of topology and combinatorics inseparable.
We conclude that one can unambiguously define a topological space $X$ of CW-type to be Buchsbaum, Cohen-Macaulay or 2-Cohen-Macaulay if $\hat{\vert}\Delta^{\wp}(X)\hat{\vert}$ is, see definitions in p.~\pageref{subsub3.3.7}.

More on these matters are found in \cite{Kod}, \cite{Mi} and \cite{W}.

\newpage
\section{Simplicial Manifolds}

\medskip
A {\bf simplicial manifold} is a {\it simplicial complex} (defined below) fulfilling certain conditions - classically involving only the local structure of the complex, usually concerning some or all ``links" or ``stars", which then are assumed to be sphere- resp. ball-like in some sense - typically topologically, homotopically, homologically or combinatorially, but there are also other ways to formulate simplicial manifold criteria.

Simplicial complexes are mainly either ``geometrical" or ``abstract". The former is usually (non-exhaustively) encountered as ``simplicial" structures, the local structure of which strongly relates to some well-behaved subspace of some Euclidian space in a way that allow a well-defined concept of ``dimension". The abstract simplicial manifolds is most likely to appear within algebraic topology or combinatorics. Through ``realization functors", ``abstract" and ``geometrical" simplicial complexes are related in some clarifying way.

Simplicial structures are useful in many ways. In particular, they allow computer-based algorithms to calculate (co)homology groups for triangulable objects, as for instance, all ``once continuously differentiable manifolds" (in particular, all ``smooth manifolds" from calculus) are triangulable, and so, any triangulation of any such compact manifold can be used to determine its (co)homology groups with respect to any homology theory due to the uniqueness theorem for simplicial homology, see \cite{E&S} p. 100.

Since this paper exclusively deals with the connection between combinatorics, general topology and algebra {\it the geometrical simplicial manifolds} will not be treated here.

\subsection{Abstract simplicial manifold}
Abstract simplicial complexes can be either ``classical" or ``augmented".
The following equality is ment to clarify the relation between ``Augmental" and ``Augmented": ``The augmental category of simplicial complexes" = ``The category of $\{\emptyset\}$-augmented simplicial complexes including both $\emptyset$ and $\{\emptyset\}$".
\begin{vardef} An abstract {\it simplicial complex} is a set S of (\underbar{non-empty}) \underbar{finite} \underbar{sets} closed under the formation of subsets.
\end{vardef}
If the parenthesis is deleted and the definition is read as {\it \underbar{non-empty} finite sets} it becomes the definition of {\it classical abstract simplicial complexes}, while when the word {\it non-empty} is completely ignored, i.e. when the definition is read with any {\it \underbar{finite sets}} -- including the empty set $\emptyset$, the definition of {\it augmented abstract simplicial complexes} is at hand.
$$\mbox{Formally: }(\sigma\in \Sigma)\land(\emptyset\neq\tau \subset \sigma)\Longrightarrow\tau\in\Sigma \ \mbox{ respectively }\ (\sigma\in\Sigma) \land(\tau\subset\sigma)\Longrightarrow\tau\in\Sigma.$$

One of the most significant difference between the classical and the augmented abstract simplicial complexes is disclosed in the definition of the binary operation of the {\it join} (${\ast}$) of two simplicial complexes. $\emptyset$ is the join-unit with respect to the classical definition, while $\{\emptyset\}$ is the join-unit with respect to the augmented definition. This difference between the join-definitions is crucial and is therefore written out together with equivalent definitions to those just given for simplicial complexes.

\medskip\noindent
{\bf Definition.} An {\it augmented abstract simplicial complex} $\Sigma$ on a {\it vertex set} $V_\Sigma$ is a collection (possibly empty)
of finite subsets $\sigma$, the {\it simplices}, of $V_\Sigma$ satisfying:
$$
\mbox{\hskip-1.5cm(a)}\mbox{ If } v\in V_{\Sigma}, \text{ then } \{v\}\in
\Sigma.
$$
$$
\mbox{(b)} \mbox{ If } \sigma\in \Sigma \text{ and }
 \tau\subset \sigma, \text{ then } \tau\in \Sigma.
$$
Define the {\it simplicial join} $\Sigma_1\ast\Sigma_2$, of two simplicial complexes $\Sigma_1$ and $\Sigma_2$ with \underbar{disjoint} vertex sets, to be:
$$
\Sigma_1\ast\Sigma_2:=\{\sigma_1\cup\sigma_2 \mid
\sigma_i\in\Sigma_i \ (i=1,2) \}.
$$

\medskip\noindent
{\bf Definition.} (See \cite{Sp} p. 108-9) A {\it classical abstract simplicial complex} $\Sigma$ on a {\it vertex set} $V_\Sigma$ is a collection of \underbar{non-empty} finite subsets $\sigma$, the {\it simplices}, of $V_\Sigma$ satisfying:
$$
\mbox{\hskip-2.26cm(a)}\mbox{ If } v\in V_{\Sigma}, \text{ then } \{v\}\in
\Sigma.
$$
$$
\mbox{(b)} \mbox{ If } \sigma\in \Sigma \text{ and }
 \emptyset\ne\tau\subset \sigma, \text{ then } \tau\in \Sigma.
$$
Define the {\it classica simplicial join} $\Sigma_1\ast\Sigma_2$, of two simplicial complexes $\Sigma_1$ and $\Sigma_2$ with \underbar{disjoint} vertex sets, to be:
$$
\Sigma_1\ast\Sigma_2:=\Sigma_1\cup\Sigma_2\cup\{\sigma_1\cup\sigma_2 \mid
\sigma_i\in\Sigma_i \ (i=1,2) \}.
$$

The importance of the augmented join can not be overrated,   since it is the functorial counterpart of the (graded) tensor produkt in commutative algebra under the Stanley-Reisner ring construction, cf. Example~{\bf iii} p. \pageref{Ex1:03}.

\subsection{Classical abstract simplicial manifolds}
$\phantom{.}$

\medskip
When combinatorialists, in the 1970th, began to use  commutative algebra to solve combinatorial problems, it became obvious that the category of {\it classical} abstract simplicial complexes were inadequate for their purposes and  they turned to the category of {\it augmented} abstract simplicial complexes, which will be the target also of the investigation below.

\subsection{Augmented abstract simplicial manifolds}
\subsubsection{Background}
"Classical abstract simplicial complexes" are present in the main-bulk of classical algebraic topological literature while the contemporary literature on simplicial complexes usually just takes the augmented concept for granted when defining or mentioning ``simplicial complexes".

The homology theory used below, named augmental homology and denoted $\hat{\mathbf{H}}_{\ast}$, is a formal homology theory that obeys the Eilenberg-Steenrod axioms and unifies the classical relative homology functor and the classical reduced homology functor. Its definition is concretizised in the chapter ``Augmental homology" and, in particular, it formalizes the informal use of the classical reduced homology within combinatorics, which sometimes is stated simply as ``We use $\tilde{\mathbf{H}}_{\ast}$ with $\tilde{\mathbf{H}}_{-1}( \{\emptyset\};\mathbf{G})=\mathbf{G}$" - a working practice that was introduced in the 70th during the algebraization of
combinatorics and which is sufficient as long as one deals with simplicial homology and is able to avoid relative homology through the use of the ``link"-construction.

The contemporary/augmented, i.e. the non-classical, definition of {\it simplicial complexes} allows two vertex-free complexes $\emptyset \mbox{ and } \{\emptyset\}$, the latter of which actually deprives modern combinatorics of a formal {\it realization functor} due to the fact that the classical category of {\it topological spaces and continuous functions} contains only one pointless space, namely $\emptyset$. This is formally rectified by performing an analogous algebraization of general topology by simply introducing a $(-1)$-dimensional topological join-unit (denoted ``$\{\wp\}$") that matches the $(-1)$-dimensional simplicial complex $\{\emptyset\}$.

There is a ``notational" issue that steams from the fact that the empty set $\emptyset$ plays a dual role in the contemporary category of simplicial complexes -- it is both a $(-1)$-dimensional \underbar{simplex} and a $(-\infty)$-dimensional \underbar{simplicial complex}. To stress those different roles,  $\emptyset$ is denoted $\emptyset_o$ when, and only when, it is used as a $(-1)$-dimensional \underbar{simplex}. The new simplicial complex $\{\emptyset_o\}$ will serve as the  ( {\it abstract simplicial}) $(-1)\text{-standard simplex}.$

\subsubsection{Definitions} Below, definitions will be given for the following three simplicial manifold structures: {\bf 1} {\it pseudomanifolds}, {\bf 2} {\it quasi-manifolds} and {\bf 3} {\it \underbar{simplicial} homology-manifolds}. The two first are combinatorial structures, though the boundary definition of quasi-manifolds isn't, as it is equivalently formulated as that for homology manifolds, which are completely homology-based and thereby dependent on the particular choice of coefficient module. The hierarchy is as follows:
$\Sigma$ is a homology manifold $\Longrightarrow$ $\Sigma$ is a quasi-manifold $\Longrightarrow$ $\Sigma$ is a pseudomanifold. None of these implications is an equivalence.\\
Also a definition {\bf 3´} for {\it \underbar{singular} homology-manifolds} is given, mainly to accentuate the connection between the combinatorial and the point-set  topological manifold structures.

\medskip\noindent
{\bf Definition 1.} An $n$-dimensional pseudomanifold is a locally finite $n$-complex
${\Sigma}$ such that;
\\
{\bf ($\alpha$)} ${\Sigma}$ is pure, i.e. the maximal simplices in ${\Sigma}$ are all $n$-dimensional.
\\
{\bf ($\beta$)} Every $(n-1)$-simplex of ${\Sigma}$ is the face of at most two $n$-simplices.
\\
{\bf ($\gamma$)} If $s$ and $s'$ are n-simplices in ${\Sigma}$, there is a finite sequence $s=s_0,s_1,\ldots s_m=s'$ of n-simplices in ${\Sigma}$ such that $s_i\cap s_{i+1}$ is an $(n-1)$-simplex for $0\le i<m$.\\
($\gamma$ is precisely the definition of {\it a pure} and  {\it strongly connected simplicial complex} if any pair of maximal simplices can be connected in this way.)

\smallskip
The boundary, $\mathrm{Bd}_{\Sigma}$, of an n-dimensional pseudomanifold ${\Sigma}$, is the subcomplex generated by those (n-1)-simplices which are faces of exactly one n-simplex in ${\Sigma}$.

\medskip\noindent
{\bf Definition 2.} \label{quasim} A quasi-$n$-manifold is a locally finite $n$-dimensional simplicial complex ${\Sigma}$ and;
\\
{\bf ($\alpha$)}${\Sigma}$ is pure. ($\alpha$ is redundant -- as a consequence of $\gamma$ (A. Björner))
\\
{\bf ($\beta$)}Every $(n-1)$-simplex of ${\Sigma}$ is the face of at most two $n$-simplices.
\\
{\bf ($\gamma$)}$\mbox{Lk}_{_{\Sigma}}\sigma$ is connected
i.e. $\hat{\mathbf{H}}_{0}({\rm Lk}_{_{\Sigma}}\sigma;\mathbf{G})=0$ for all
$\sigma\in\Sigma,$ s.a. $\dim\sigma<n-1,$ i.e. $\dim{\mathrm{Lk}}_{_{^{\Sigma}}}\sigma\geq1.$

\smallskip
The boundary with respect to the coefficient module {\bf G} of a quasi-{\it n}-manifold $\Sigma$, denoted $\mbox{Bd}_{_{{\bf G}}}\Sigma,$ is the complex;
$$
\mathrm{Bd}_{_{{\bf G}}}\Sigma:= \{\sigma\in\Sigma\mid\hat{\mathbf{H}}_{n-\#\sigma}
(\mathrm{Lk}_{\Sigma}\sigma;\mathbf{G})=0\} \ (\cong\{\sigma\in\Sigma\mid\hat{\mathbf{H}}_{n}(\Sigma,\mathrm{cost}_{\Sigma}\sigma;\mathbf{G})=0\}),
$$
where {\bf G} is a unital module over a commutative ring ${\mathbf{A}}.$ ``Lk" and ``cost" are defined
in p.~\pageref{DefLk}.
$\mathrm{Bd}_{_{{\bf G}}}\Sigma \mbox{ is } \mathbf{G}$-dependent and therefore, in general, \underbar{not} purely combinatorial.

\medskip\noindent
{\bf Lemma.} $\Sigma$ is a quasi-$n$-manifold if all of its links $\mathrm{Lk}_{\Sigma}\sigma$ are pseudomanifolds. (In particular including $\mathrm{Lk}_{\Sigma}\emptyset_o=\Sigma.$)

Recall that $\underline{\dim\mbox{Lk}_{\Sigma}\sigma=\dim\Sigma-\#\sigma}$. 
In the next definition, let: $n:=\dim\Sigma.$

\medskip\noindent
{\bf Definition 3.} $\emptyset$ is defined to be a (simplicial) $\mbox{homology}_{_{{\bf G}}}$-$(-\infty)$-manifold and $\Sigma={\bullet}{\bullet}$ is a $\mbox{homology}_{_{{\bf G}}}0$-manifold.
Any other (than $\emptyset \mbox{ or } {\bullet}{\bullet}$) connected, locally finite $n$-complex in ${\mathcal{K}}_{\wp}$ is a homology$_{_{{\bf G}}} n\mbox{-manifold }$ $(n\mbox{-hm}_{_{{\bf G}}})$ if
$$\begin{array}{ll}
\hat{\mathbf{H}}_{i-\#\sigma}
(\mathrm{Lk}_{\Sigma}\sigma;\mathbf{G})=0  \mbox{ if } i\ne n \mbox{ for all } \emptyset_o\ne \sigma\in \Sigma,
&(\mathbf{i})\\
\hat{\mathbf{H}}_{n-\#\sigma}(\mathrm{Lk}_{\Sigma}\sigma;\mathbf{G})
\cong 0 \mbox{ or }\mathbf{G} \mbox{ for all } \emptyset_o\ne \sigma\in \Sigma.
&(\mathbf{ii})\end{array} \hskip3cm(\mathbf{1})$$
The boundary is defined to be:
$\mbox{Bd}_{_{{\bf G}}} \Sigma:= \{\sigma\in \Sigma \mid \hat{\mathbf{H}}_{n-\#\sigma}(\mathrm{Lk}_{\Sigma}\sigma;{\bf G})= 0\}.$

\noindent
$\mbox{An }n\mbox{-hm}_{_{{\bf G}}} \Sigma$ is {\it joinable} if $(\mathbf{i})$ holds also for $\sigma=\emptyset_o$.

\noindent
$ \mbox{An }n\mbox{-hm}_{_{{\bf G}}} \Sigma\not=\emptyset
\mbox{ is a homology}_{_{{\bf G}}}\mbox{-}n\mbox{-} \mbox{ sphere }
(n{\mbox{-hsp}}_{_{{\bf G}}})$ $ \mbox{ if, for all }\sigma\in \Sigma,\mbox{ } \hat{\mathbf{H}}_{n-\#\sigma}(\mathrm{Lk}_{\Sigma}\sigma;\mathbf{G})
={\mathbf{G}}$.

\medskip\noindent
{\bf Definition 4.} Let ``manifold" stand for a pseudo-, quasi- or a homology manifold.
A finite/compact $n$-manifold ${\mathcal{S}} \mbox{ is } \mbox{orientable}_{_{{\bf G}}}$ if $\hat{\mathbf{H}}_{n}({\mathcal{S}},{\mathrm{Bd}}{\mathcal{S}};{\mathbf{G}})\cong{\mathbf{G}}.$
An n-manifold is $\mbox{orientable}_{_{{\bf G}}}$ if all its finite/compact $n$-submanifolds are orientable $-$ else, non-orientable$_{_{{\bf G}}}$.
Orientability is left undefined for $\emptyset$.

\smallskip
To be able to analyze the boundary of a simplicial manifold in detail one needs to be able to identify {\it the strongly connected boundary components}, see Definition 1$\gamma$.

\medskip\noindent
{\bf Definition 5.} $\{{\bf B}^{\Sigma}_{{\bf G},j}\}_{j\in{\bf I}}$ is {\it the set of strongly connected boundary components of}  $\Sigma$ if
$\{\mathbf{B}^{\Sigma}_{{\bf G},j}\}_{j\in{\bf I}}$ is the maximal strongly connected components of ${\rm Bd}_{{\bf Bd}_{{\bf G}}}\Sigma$.
$(\Rightarrow\mathbf{B}^{\Sigma}_{{\bf G},j}$ is pure and if $\sigma$ is a maximal simplex in $\mathbf{B}_{j}\ (:=\mathbf{B}^{\Sigma}_{{\bf G},j})$, then;
$\mathrm{Lk}_{({\rm Bd}_{{\bf G}}\Sigma)}\sigma
=\mathrm{Lk}_{{\bf B}_j} \sigma=\{\emptyset_o\}.)$

\medskip\noindent
{\bf Note.} The functor
$$ {\mathcal{F}}_{\wp}:\ {\mathcal{D}}\longrightarrow {\mathcal{D}}_{\wp};\ {\mathcal{F}}_{\wp} (X)=X+\{\wp\}=:X_{\wp}$$
simply adds, using ``topological sum", the $(-1)$-{\it dimensional topological join-unit} $\{\wp\}$ to each classical topological space in $ {\mathcal{D}} $ and it is defined in the chapter on Augmental homology.
$$ {\mathcal{F}}_{\wp} \mbox{ has an invers, } {\mathcal{F}}:\ {\mathcal D}_{\wp}:\longrightarrow {\mathcal D};\ {\mathcal{F}}(X_{\wp})=X \mbox{ deleting }{\wp}.$$
Note also that ${\mathcal{F}}(\{\wp\}) = {\mathcal{F}}(\emptyset) = \emptyset.$

So, the typical augmental topological space is a construction of the following kind: $X_{\wp}={\mathcal F}_{\wp}(X)=X+{\wp}$, where ``$X$" is a classical topological space and the ``$+$" denotes ``{\it topological sum}".
To be able to discuss general topological questions concerning augmental topological spaces the following general doctrine is formulated as a guideline:

\medskip\noindent
{\bf Doctrine:} The underlying principle for definitions is that a concept in
${\mathcal D}\ ({\mathcal{K}})$ is carried over to ${\mathcal{D}}_{\wp}\ ({\mathcal{K}}_o)$ by ${\mathcal{F}}_{\wp}\ ({\mathcal{E}}_o)$ with addition of explicit definitions of the $\mbox{concept}_{o}$ for cases which are not proper images under
$ {\mathcal{F}}_{\wp}\ ({\mathcal{E}}_o)$.

\smallskip
Whatever manifold $\mathbf{M}$ is at hand, its boundary in the augmented setting is essentially the same as in the classical setting since;

\smallskip
${{{\mathcal{E}}}}({{\rm Bd}}\mathbf{M})={{\rm Bd}}({{{\mathcal{E}}}}(\mathbf{M}))$ in the simplicial case,

\noindent
resp.

${{{\mathcal{F}}}}({{\rm Bd}}\mathbf{M})={{\rm Bd}}({{{\mathcal{F}}}}(\mathbf{M}))$ in the topological case.

\medskip
So, what remains after removing $\emptyset_o$ resp. $\wp$ from a boundary in the augmented setting is exactly the classical boundary - always!

In the opposite direction: For any classical $\mbox{homology}_{_{{\bf Z}}}$ manifold
$X\neq\bullet$ (:= The  one-point space) with $\mbox{Bd}_{_{{\bf Z}}}X=\emptyset$;
Bd$_{_{{\bf Z}}}{\mathcal F}_{\wp}(X) =\emptyset$ if
$X$ is compact and orientable and Bd$_{_{{\bf Z}}}{\mathcal{F}}_{\wp}(X) = \{\wp\} \mbox{ else}$.

The one-point space (resp. the one-vertex complex) $\bullet$ is the only compact \underbar{orientable} manifold with boundary $=\{\wp\}\ (\{\emptyset_{o}\}).$

The subindex $\wp$ will usually be deleted, as for instance, $X_\wp$ will usually be written $X.$

\subsubsection{Realizations}A topological space $X$ that is homeomorphic to the realization $|\Sigma|$ of an abstract simplicial complex $\Sigma$ is called a {\it polytope} or a {\it triangulable space} and $\Sigma$ is said to be a {\it triangulation} of $X$. The augmented realization functor essentially equals the classical, except that it also maps $\{\emptyset_o\} \mbox{ onto } \{\wp\}$, cf. \cite{Fo1} p.~13.\\
A topological space is called a $n$-{\it pseudomanifold} or a {\it quasi-$n$-manifold} if it can be triangulated into a abstract simplicial complex that is a $n$-{\it pseudomanifold} resp. a {\it quasi-$n$-manifold} as defined through the above two purely combinatorial structures. It can be shown that if $X$ has a triangulation $\Sigma$ that is a pseudo-/quasi-$n$-manifold then any triangulation of $X$ is a pseudomanifold or a quasi-$n$-manifold resp.

A different definition-approach is used when declaring a topological space to be a (singular) $\mbox{homology}_{_{{\bf G}}} n\mbox{-manifold }$ $(n\mbox{-hm}_{_{{\bf G}}})$.

Recall the definition of the
``extended setminus" $\setminus_o$ p.~\pageref{DefSetminus} and that: $\Sigma$ is (locally) finite $\Longleftrightarrow \vert\Sigma\vert$ is (locally) compact.

\medskip\noindent
{\bf Definition 3´.} $\emptyset$ is said to be a (singular) $\mbox{homology}_{_{{\bf G}}}\mbox{-}(-\infty)$-manifold and $X={\bullet}{\bullet}$ is a $\mbox{homology}_{_{{\bf G}}}0$-manifold.
Any other (than $\emptyset$ or ${\bullet}{\bullet}$) connected, locally compact Hausdorff space $X\in{{\mathcal D}}_{\wp}$ is a (singular)
homology$_{_{{\bf G}}} n$-manifold
(n-hm$_{_{{\bf G}}}$) if
$$\begin{array}{ll}
\hat{\mathbf{H}}_{i}(X,X\setminus_{o}x;\mathbf{G})
 = 0  \mbox{ if } i\ne n \mbox{ for all } \wp\ne x\in X, & (\mathbf{i})\\
\hat{\mathbf{H}}_{n}(X,X\setminus_{o}x;{\mathbf{G}})
\cong 0 \mbox{ or }\mathbf{G}\mbox{ for all } \wp\ne x\in X
\mbox{ and} = \mathbf{G}\mbox{ for some } x\in X. & (\mathbf{ii})
\end{array} \hskip1.5cm({\mathbf{1}}^{\prime})$$
The boundary is defined to be:
$$\mbox{Bd}_{_{{\bf G}}} X:= \{x\in X\mid{\hat{{\bf H}}}_{n}(X,X\setminus_{o} x;{\mathbf{G}})=0\}.$$
$\mbox{An }n\mbox{-hm}_{_{{\bf G}}} X$ is {\it joinable} if $(\mathbf{i})$  holds also for $x=\wp.$
$$\hskip-0.2cm\mbox{An }n\mbox{-hm}_{_{{\bf G}}} X\not=\emptyset
\mbox{ is a homology}_{_{{\bf G}}}\mbox{-}n\text{-sphere} \  (n\mbox{-hsp}_{{\bf G}}) \mbox{ if, for all }x\in X, \hat{{\bf H}}_{n}(X,X\setminus_{o}x;{\mathbf{G}})={\mathbf{G}}.$$

$\mbox{A simplicial complex } \Sigma\mbox{ is a hm}_{_{{\bf G}}} \mbox{ if } |\Sigma| \mbox{ is}.$
So, in particular, a triangulable $n\mbox{-hsp}_{_{{\bf G}}}$ is a compact space.

Note that even if $X\mbox{ is a hm}_{_{{\bf G}}}$ it may not be triangulable.

The following proposition implies that the boundary definition for quasi- and homology manifolds resp. are identical for polytopes.

The above mentioned $(-1)$-dimensional topological join-unit, ``$\wp$", shows up in different disguises, e.g., in the next proposition it is incarnated as the ``{\it trivial barycentric vertex-function}" $\alpha_{0}$ fulfilling $\alpha_{0}({v})\equiv0$ for all vertices ${v}$.

Since ``$\mbox{Lk}_{_{\Sigma}}\emptyset=\Sigma$", the use of this new unifying {\it relative homology} functor together with the  ``natural" (in particular faithful)
{\it realization functor} now requires a minor modification/extension p.~\pageref{DefSetminus},
denoted ``$\setminus_{o}$", of the classical {\it setminus} ``$\setminus$" resulting in the following very useful proposition. (See the chapter on Augmental homology.)

\begin{varprop} \label{Prop:00} Let $\mathbf{G}$ be a $($unital$)$ module  over a commutative ring {\bf A} with unit.
With $\alpha$ in the interior of $\sigma$, i.e. $\alpha\in\{\beta\in|\Sigma|\mid[\mbox{v}\in\sigma] \Longleftrightarrow [\beta(\mbox{v})\neq0]\}$ and $\alpha=\alpha_{0}$ if and only if
$\sigma=\emptyset_o$,
the following {\bf A}-module isomorphisms are all induced by chain $($homotopy$)$ equivalences.
$$\hat{\mathbf{H}}_{i-\#\sigma}(\mathrm{ Lk}_{\Sigma}\sigma;\mathbf{G})\cong
\hat{\mathbf{H}}_{i}(\Sigma,\mathrm{ cost}_{\Sigma}\sigma;\mathbf{G})\cong\hat{\mathbf{ H}}_{i}(\vert\Sigma\vert,\vert\mathrm{ cost}_{\Sigma}\sigma\vert;
\mathbf{G})\cong\hat{\mathbf{ H}}_{i}(|\Sigma|,|\Sigma|\setminus{_o}\alpha;\mathbf{G}).$$
\end{varprop}
\subsubsection{Product and Joins of manifolds}Three types of manifolds, complexes \& polytopes and two binary operation makes 12 different situations that below are rationalized into one theorem. Some preparation is needed though.
The (ordered simplicial) Cartesian product of simplicial complexes are given at the end of the chapter on Stanley-Reisner rings, while the join-definition is given in the beginning of that article.
\\
When $\diamond$ in Theorem 01 below is interpreted throughout as product $\times$, the word manifold temporarily excludes $\emptyset,\{\emptyset_o\}$ and also the $0$-sphere $\bullet\bullet$, since the latter in general makes the resulting product non-connected. $\emptyset$ is excluded for joins $\ast$. For products it is assumed that $\epsilon:=0$, while $\epsilon:=1$ in the case of joins.
For joins, let the word manifold on the right hand side of Th. 01.1 be limited to ``any \underbar{finite} manifold" respectively to ``any compact \underbar{joinable} homology$_{_{{\bf G}}} n_{i}$-manifold" when homology manifolds are concerned.

Local compactness is preserved under join if, and only if, the factors are both compact (and regular according to D.E. Cohen, 1956).

For quasi- and homology manifolds, ``Bd" should be interpreted as ``$\mathrm{Bd}_{_{{\bf k}}}$", with $\mathbf{k}$ a field, in the following theorem.

\begin{theorem}\label{Th01P22}
$\mbox{For polytopial spaces } X_{1}, X_{2} \mbox{ and any field }\mathbf{k}:$
$$\hskip-2.1cm\mathbf{1.}\ X_{1}\diamond X_{2}\ \text{is a}\ (n_{1}+n_{2}+\epsilon)
\mbox{-manifold}\Longleftrightarrow X_{i}\mbox{ is a }n_{i}\mbox{-manifold for } i=1,2.$$
$$\mathbf{2.}\ {\rm Bd}(\bullet\times X)= \bullet\times({\rm Bd} X).\mbox{ Otherwise } {\rm Bd}(X_{1}\diamond X_{2})=(({\rm Bd}
X_{1})\ \diamond X_{2})\cup(X_{1}\diamond({\rm Bd}X_{2})).$$
$$\hskip-5.0cm\mathbf{3.}\ X_{1}\diamond X_{2}\mbox{ is orientable}_{_{{\bf k}}}
\Longleftrightarrow\mbox{ both } X_{1},X_{2}\mbox{ are orientable}_{_{{\bf k}}}.$$
\end{theorem}
Theorem~\ref{Th01P22} can be generalized considerably through Theorem~\ref{Th06} p.~\pageref{Th06} in the chapter on Augmental Homology.
$$\mbox{\bf Note. } (X_{1}\diamond X_{2}, {\rm Bd}_{k}(X_{1} \diamond X_{2}))=
\begin{array}{l}\left[{\text{\rm Motivation: Th.~\ref{Th01P22}.2 +}\atop{\text{the pair-}\times{\tiny/}\ast\text{-formulas}}}\right]\end{array}
=(X_{1},{\rm Bd}_{_{{\bf k}}}{}_{}X_{1})
\diamond
(X_{2},{\rm Bd}_{_{{\bf k}}}X_{2}).
$$

\subsubsection{Examples}
$\phantom{.}$

{\bf Example 1.} There are a number of homeomorphisms involving products and joins of {\it m}-unit balls/disks $\mathbb{E}^m$ and $m$-unit spheres $\mathbb{S}^m$. For example, with $n:=p+q, \ p,q\ge0$ then,
$$\mathbb{E}^{p}\ast \mathbb{E}^{q}\simeq\mathbb{E}^{n+1}\simeq\mathbb{E}^{p}\ast\mathbb{S}^{q},\indent \indent
\mathbb{S}^{n+1}\simeq\mathbb{S}^{p}\ast \mathbb{S}^{q} \indent {\rm and}\indent \mathbb{E}^{p}\times \mathbb{E}^{q}\simeq\mathbb{E}^{n}.$$
Now the above boundary formula provides some of the homeomorphisms used within surgery theory, as for instance:
$$\mathbb{S}^{_{^n}}=\mbox{Bd}\mathbb{E}^{n+1}\simeq\mbox{Bd}
(\mathbb{E}^{p}\ast\mathbb{E}^{q})\simeq\mathbb{E}^{p}\ast
\mathbb{S}^{q-1}\cup\mathbb{S}^{p-1}\ast\mathbb{E}^{q}\simeq
\mbox{Bd}(\mathbb{E}^{p+1}\times\mathbb{E}^{q})
\simeq\mathbb{E}^{p+1}\times\mathbb{S}^{q-1}\cup \mathbb{S}^{p}\times \mathbb{E}^{q}.$$

\smallskip\noindent 
{\bf Example 2.} The real projective plane is the boundary of the cone of a Möbius band $\mathbb{M}$ when regarded as a quasi-manifold (or as a pseudomanifold), and using, for instance, the prime field $\mathbb{Z}_3$ with three elements as coefficient group.
$(\mathbb{M}$ is not a {\it joinable} 2-hm$_{{\bf k}}.$ $\mathbb{RP}^2$  is a {\it joinable} 2-hm$_{{\bf k}}\text{ if }\mbox{ char}{\mathbf{k}}\ne2$. See \cite{Gr} p.~36.)
$$\mbox{Bd}_{\mathbb{Z}_3}(\mathbb{M}\ast\bullet)=
(({\rm Bd}_{\mathbb{Z}_3} \mathbb{M} )\ast\bullet)\cup(\mathbb{M}\ast{\rm Bd}_{\mathbb{Z}_3}\bullet)
=(\mathbb{S}^1\ast\bullet)\cup(\mathbb{M}\ast{\rm Bd}_{\mathbb{Z}_3}\bullet)=
{\mathbb{E}}^2 \cup(\mathbb{M}\ast \{\emptyset_o\})={\mathbb{E}}^2 \cup{\mathbb{M}},$$
$\mbox{which is a well-known century-old representation of the real projective plane } \mathbb{RP}^2$.

\medskip\noindent
{\bf Note.} $\emptyset, \{\emptyset_o\},$ and $0$-dimensional complexes
with either one, ${\bullet},$ or two, ${\bullet}{\bullet},$ vertices
are the only manifolds in dimensions $\le0$, and the
$|1$-manifolds$|$ are finite/infinite $1$-circles and (half)lines, while,
$$[{\Sigma} \mbox{ is a quasi-2-manifold}] \Longleftrightarrow
[\Sigma \mbox{ is a homology}_{\mathbb{Z}} \mbox{ 2-manifold}].$$
Any augmental topological space $X_{\wp}=X+{\wp}$ is by construction, as a topological sum, a non-connected topological space, but through the above accepted doctrine it is connected if and only if $X$ is.
For example -- as stated above -- Definition $\mathbf{ 1.\gamma}$ is paraphrased as  ``$\Sigma$ is {\it strongly connected}".
Now, ${\bullet}{\bullet}$ is strongly connected as a simplicial complex, due to existence of the $(-1)$-dimensional simplex $\emptyset_o$, while its realization $\vert{\bullet}{\bullet}\vert$, usually also denoted simply ${\bullet}{\bullet}$, is non-connected, due to the above doctrine.

Note also that
$\mathbb{S}^{-1} :=\{\emptyset_o\}$ is the boundary
of the $0$-ball $\bullet$, the double of which is the $0$-sphere $\bullet\bullet$. Both the $(-1)$-sphere
$\{\emptyset_o\}$ and the $0$-sphere $\bullet\bullet$ has, as preferred, empty boundaries.

\subsubsection{Weak manifolds in general topology}
The general topological concept of a ``{\it weak homology manifold}" is essentially equivalent to that of ``{\it Buchsbaum complexes}" within combinatorics. This illustrates the usefulness of the above mentioned algebraization of general topology. Of course, this connection goes both ways and so, further motivates the contemporary effort to topologize algebraic structures.
\\
The general topological regularity condition ``Hausdorff" in the following definition of a ``{\it weak homology manifold}" is equivalent to $\mbox{T}_2 + \mbox{T}_1$, see for instance Wikipedia: ``Separation axioms".

\medskip\noindent
{\bf Definition 6´.} $\emptyset$ is defined to be a ${\it weak}\ {\it homology}_{_{{\bf G}}}$-$(-\infty)$-manifold.
A nonempty locally compact Hausdorff space
$X\in{\mathcal D}_{\wp}$ is a {\it weak homology}$_{_{{\bf G}}} n\text{-manifold }(n$-whm$_{_{{\bf G}}})$ if, for some $\mathbf{A}$-module ${\mathcal{R}}$;
$$\begin{array}{ll}
\hat{\mathbf{H}}_{i}(X,X\setminus_{o}x;\mathbf{G}) = 0 \text { if } i\ne n \text{ for all }
\wp\ne x\in X,
& (\mathbf{i})\\
\hat{\mathbf{H}}_{n}(X,X\setminus_{o}x;\mathbf{G})
\cong
\mathbf{G}\oplus {\mathcal{R}} \text{ for some } \wp\ne x\in X \text{ if } X\ne\{\wp\}.
& (\mathbf{ii}^\prime)
\end{array}\hskip2cm({\mathbf{2}}^{\prime})$$
$\hskip-0.2cm\begin{array}{ll}
\text{An }n\mbox{-whm}_{_{\!{\bf G}}}X\text{ is joinable } (n\mbox{-jwhm}_{_{\!{\bf G}}})\text{ if } (\mathbf{i}) \text{ holds also for } x=\wp. \\
\text{An }n\mbox{-jwhm}_{_{\!{\bf G}}} X
\text{ is a weak homology-}n\mbox{-sphere}_{_{\!{\bf G}}} (n\mbox{-whsp}) \text{ if }
\hat{\mbox{H}}_{n-1}(X\!\setminus_{o}x;\mathbf{G})=0 \text{ for all }x\!\in\! X.
\end{array}$

\medskip\noindent
{\bf Definition.} $
X \mbox{ is {\it acyclic}}_{_{\!{\bf G}}} \mbox{if }
\hat{\mathbf{H}}_{i}(X,\emptyset;{\mathbf{G}})=0 \mbox{ for all } i \in {\mathbb{Z}}.\
\text{ So, } \{\wp\} (=|\{\emptyset_o\}|) \text{ isn't acyclic}_{_{\!{\bf G}}}\!.
$

\noindent
$\text{An arbitrary }X\mbox{ is {\it ordinary}}_{_{\!{\bf G}}}\mbox{if }{\mathbf{2}}^{\prime}\text{ implies that }\hat{\mathbf{H}}_{i}(X\setminus_{o}x;\mathbf{G})=0, \mbox{ for all }i\ge n$\mbox{ and all $x\!\in\! X.$}

$X \mbox{ is {\it locally weakly direct}}_{{\bf G}} \mbox{ if }
\hat{\mathbf{H}}_{i}(X,X\setminus_{o}x;{\mathbf{G}})
\cong
\mathbf{G}\oplus {\mathcal{Q}}
\mbox{ for some } i, \mbox{ some } {\mathbf{A}}\mbox{-module } {\mathcal{Q}}$
and some $\wp\ne x\in X $.
(A definition of technical nature - used to avoid complete annihilation of certain graded tensor products.)

$\hat{\mathbf{H}}_{_{\dim{\Sigma}}} (|\Sigma|,|\Sigma|\setminus_{o}\alpha;\mathbf{G})\cong{\mathbf{G}}$
if $\alpha\in \mbox{Int}\sigma$ and $\sigma$
is a simplex of maxidimension i.e., if $\#\sigma-1=:\dim\sigma=\dim\Sigma$, since now
Lk$_{_{\Sigma}}\sigma=\{\emptyset_{_{^{o}}}\}$.
So, for polytopes, ($\mathbf{ii}^\prime$) is always fulfilled,
which simplifies the definition of {\it weak homology manifolds} in the category of simplicial complexes, since only condition ($\mathbf{i}$) needs to be justified i.e.:

\medskip\noindent
{\bf Definition 6.} $\emptyset$ is defined to be a ${\it weak}\ {\it homology}_{_{{\bf G}}}$-$(-\infty)$-manifold.\\
A nonempty simplicial complex $\Sigma\in{\mathcal{K}}_{o}$ is a {\it weak homology}$_{_{{\bf G}}} n\text{-manifold } (n$-whm$_{_{{\bf G}}})$ if;
$$
\hat{\mathbf{H}}_{i-\#\sigma}
(\mathrm{Lk}_{\Sigma}\sigma;\mathbf{G})=0  \mbox{ if } i\ne n \mbox{ for all } \emptyset_o\ne \sigma\in \Sigma,
\hskip3cm(\mathbf{i})\hskip2.0cm(\mathbf{2})$$
$\begin{array}{ll}
\text{An }n\text{-whm}_{_{{\bf G}}} \Sigma \text{ is joinable }
(n\text{-jwhm}_{_{{\bf G}}}) \text{ if } (\mathbf{i}) (\mathbf{2}) \text{ holds also for } \sigma=\emptyset_o. \\
\text{An } n\text{-jwhm}_{_{{\bf G}}} \Sigma  \text{ is a weak homology-}n\text{-sphere}_{_{{\bf G}}} (n\text{-whsp}) \text{ if } \hat{\mbox{H}}_{n-1}(\mathrm{cost}_{\Sigma}\sigma;\mathbf{G})=0 \text{ for all } \sigma\in \Sigma.
\end{array}$

\medskip\noindent
{\bf Note.} Nonempty triangulable manifolds are {\it ordinary} if the coefficient module is a field$\mbox{ or }{\mathbb{Z}}.$

\begin{theorem}\label{Th02P24} $\mbox{For ordinary}_{{\bf k}}$ {\it locally weakly direct}$_{{\bf k}}
$
${\mathbf{T}}_1$-spaces $X_{1},X_{2}\ (X_{i}\neq \emptyset, \{\wp\});$

\smallskip\noindent
{\bf i.} $X_{1}\times X_{2} $
$(n_{1}+n_{2})\mbox{-whm}_{{\bf k}}$
$\Longleftrightarrow X_{1},X_{2}\mbox{ both whm}_{{\bf k}}.
$

\smallskip\noindent
{\bf ii.} If $n_{_{^1}} +n_{_{^2}}>n _{_{^i}}\ i=1,2$,
then,
$$[X_{1}\times X_{2}
(n_{1}+n_{2})\mbox{-jwhm}_{{\bf k}}]\!\Longleftrightarrow [X_{1},X_{2}
\mbox{ both } n_{_{^i}}\mbox{-jwhm}_{{\bf k}}\mbox{ and acyclic}_{{\bf k}}].
$$
$(\mbox{Th.~\ref{Th02} p.~\pageref{Th02} in ch.~2 on augmental homology implies that;}$
$$[\hat{\mathbf{H}}_{_i}(X_{_{1}}{\times}X_{_{2}};{\mathbf{k}})=0
\mbox{ for }i\not=n_{1}+n_{2}]\Longleftrightarrow
[X_{_{1}}, X_{_{2}}\mbox{ both acyclic}_{{\bf k}}]\Longleftrightarrow
[X_{_{1}}{\times}X_{_{2}}\mbox{ acyclic}_{{\bf k}}].$$
$ \mbox{So, }{X_{_{1}}\times}X_{_{2}}\mbox{ is \underline{never} a whsp}_{{\bf k}}).$

\medskip\noindent
{\bf iii.} $
X_{1}\ast X_{2}
(n_{1}+n_{2}+1)\mathrm{-whm}_{{\bf k}}
$
$ \Longleftrightarrow
X_{1},X_{2} \mbox{ both }n_{_{{i}}}\mbox{-jwhm}_{{\bf k}}
$
$
\Longleftrightarrow X_{1}\ast X_{2}
\mbox{ jwhm}_{{\bf k}}.
$

\smallskip\noindent
{\bf iv.} $X_{1},X_{2}\mbox{ are both whsp}_{{\bf k}}
\mbox{ if and only if } X_{1}\ast X_{2} \mbox{ is a whsp}.$
\end{theorem}

\subsubsection{Weak manifolds in combinatorics}\label{subsub3.3.7}Combinatorialists call a finite simplicial complex a Buchsbaum ($\mathrm{Bbm}_{_{{\bf k}}}$) a  Cohen-Macaulay ($\mathrm{CM}_{_{{\bf k}}}$), a 2-Cohen-Macaulay ($2\mbox{-CM}_{_{{\bf k}}}$) or a Gorenstein (${\rm Gor}_{_{{\bf k}}}$) complex if its Stanley-Reisner ring is a Buchsbaum, Cohen-Macaulay, $2$-Cohen-Macaulay or Gorenstein ring resp. The subindex $\mathbf{k}$ indicate the coefficient module/base ring, which is limited to be a field or $\mathbb{Z}$ (for algebra-theoretical reasons).
Unlike Gorensteinness, the property of being either Buchsbaum, Cohen-Macaulay or 2-Cohen-Macaulay is a topological property in the sense that they are not sensitive to any particular triangulation of the polytope (i.e. of the realisation).
The hierarchy among these algebra founded non-triangular sensitive concepts looks as follows: $2\mbox{-CM}_{_{{\bf k}}}\Longrightarrow$CM$_{_{{\bf k}}}\Longrightarrow$ ${\rm Bbm}_{_{{\bf k}}}.$
\\
Not only polytopes, but every topological spaces is individually related to a well-defined Stanley-Reisner ring through the triangulability of the Milnor realization $\hat{\vert}\Delta^{\wp}(X)\hat{\vert}$ of the (augmented) singular complex with respect to $X$, see \cite{Mi}, as it, with respect to ${\rm Bbm}_{_{{\bf k}}}$, ${\rm CM}_{_{{\bf k}}}$ and $2\mbox{-}{\rm CM}_{_{{\bf k}}}$-ness, is triangulation invariant.
\\
Each of these four properties of simplicial complexes has been given equivalent definitions formulated purely in terms of certain homology groups of the simplicial complex at hand.
The failure of Gorensteinness to be a topological property depends on the fact that here the homology-calculations are performed on a subcomplex called {\it the core}, the definition of which depends on the {\it cone points}, which in turn depend on the triangulation of the original complex. A vertex is a ``cone point" in $\Sigma$ if it belongs to every maximal simplex of $\Sigma$.
$$\mbox{ The core of } \Sigma:=\mbox{core}\Sigma:=\{\sigma\in\Sigma \mid \sigma \mbox{ contains no cone points}\}.$$
The following three propositions, found in \cite{Sta} together with proof-references, can be used as ``definitions" when restricted to finite simplicial complexes. The proof of Proposition 3.iv is given by A. Björner as an appendix in \cite{Fo2}.
\\
Recall the definition of the $\mathcal{C}_o$-extended setminus ${\setminus_o}$ and note that the ``$i<\dim\cdot$" are all equivalent to ``$ i\ne\dim\cdot$", since there are no simplices above dimension ``$\dim\cdot$".

\begin{varprop}~{\bf 1.}\label{Prop:01} {\rm(Schenzel) (\cite{Sta} Th. 8.1)} Let $\Sigma$ be a finite simplicial complex and let $\mathbf{k}$ be a field. Then the following are equivalent:\\
$\begin{array}{llll}
\hskip-0.0cm{\bf(i.)}\ \Sigma \mbox{ is Buchsbaum over } \mathbf{k}.\\
\hskip-0.2cm\begin{array}{ll}{\bf(ii.)}\ \Sigma \mbox{ is pure, and } \mathbf{k}\langle\Sigma\rangle_{\mathfrak{p}}
\mbox{ is Cohen-Macaulay for all prime ideals ${\mathfrak{p}}$ different from}\\ \ \  \mbox{the unique homogeneous maximal ideal i.e., the irrelevant ideal.}\end{array}\\
{\bf (iii.)}\ \mbox{ For all } \sigma\in \Sigma,\ \sigma\ne \emptyset \mbox{ and } i< \dim({\rm Lk}_{\Sigma}\sigma), {\hat{\mathbf{H}}}_{i}({\rm Lk}_{\Sigma}\sigma;\mathbf{k})=0.\\
{\bf (iv.)}\ \mbox{ For all } {\alpha}\in|\Sigma|,\ \alpha\ne \alpha_0 \mbox{ and } i<\dim\Sigma,{\hat{\mathbf{H}}}_{i}(|\Sigma|,|\Sigma|\setminus_o\alpha ;\mathbf{k})=0.
\end{array}$
\end{varprop}

\begin{varprop}~{\bf 2.}\label{Prop:02} {\rm(\cite{Sta} Prop. 4.3, Cor. 4.2)} Let $\Sigma$ be a finite simplicial complex and let $\mathbf{k}$ be a field. Then the following are equivalent:\\
$\begin{array}{lll}
\hskip-0.0cm{\bf (i.)}\hskip2.6cm \Sigma \mbox{ is Cohen-Macaulay over } \mathbf{k}.\\
\hskip-0.0cm{\bf (ii.)}\ \ {\rm (Reisner)}\hskip0.6cm \mbox{ For all } \sigma\in \Sigma \mbox{ and } i< \dim({\rm Lk}_{\Sigma}\sigma), {\hat{\mathbf{H}}}_{i}({\rm Lk}_{\Sigma}\sigma;\mathbf{k})=0.\\
\hskip-0.0cm{\bf (iii.)}\ ({\rm Munkres})\indent \mbox{ For all } {\alpha}\in|\Sigma| \mbox{ and } i<\dim\Sigma,
{\hat{\mathbf{H}}}_{i}(|\Sigma|,|\Sigma|\setminus_o\alpha;\mathbf{k})=0.
\end{array}$
\end{varprop}

\begin{varprop}~{\bf 3.}\label{Prop:03}
{\rm(\cite{Sta} Th. 5.1)} Let $\Sigma$ be a finite simplicial complex, $\mathbf{k}$ a field or $\mathbb{Z}$ and $\Gamma:=\mbox{core}\Sigma.$  Then the following are equivalent:\\
\noindent
$\begin{array}{lllll}
\ {\bf (i)}\  \Sigma \mbox{ is Gorenstein over } \mathbf{k}.\\
{\bf (ii)}\mbox{ For all }  {\sigma}\in\Gamma,\ {\hat{\mathbf{H}}}_{i}(\mathrm{Lk}_\Gamma\sigma;\mathbf{k})\cong
{\Big\{}\begin{array}{ll}
\mathbf{k}\ \ \mbox{if}\ i=\dim({\rm Lk}_{\Gamma}\sigma),\\
\mathbf{0}\ \ \mbox{if}\ i<\dim({\rm Lk}_{\Gamma}\sigma)
\end{array}\\
{\bf (iii)}\mbox{ For all }  \alpha\in|\Gamma|,\ {\hat{\mathbf{H}}}_{i}(|\Gamma|,|\Gamma|\setminus_o \alpha;\mathbf{k})={\Big\{}
\begin{array}{ll}
\mathbf{k}\ \ \mbox{if}\ i=\dim\Gamma,\\
\mathbf{0}\ \ \mbox{if}\ i<\dim\Gamma.
\end{array}\\
{\bf (iv)}\ \mbox{\rm(A. Björner)}\ \Sigma \mbox{ is Gorenstein}_{\bf k}  \Longleftrightarrow\Sigma \mbox{ is C-M}_{\bf k}, \mbox{ and } \Gamma\mbox{ is an orientable pseudomanifold}\\ \mbox{without boundary.}\\
 {\bf(v)}\ { Either}\
\hskip-0.0cm\left\{\begin{array}{ll}
{\bf 1.}\ \Sigma=\{\emptyset\}, \bullet, \bullet\bullet, or,\\
\vbox{\hsize13cm\hskip-0.4cm{\bf 2.}\ $\Sigma$ is Cohen-Macaulay over $\mathbf{k}$,  $\dim\Sigma\geq 1$, and the link of every $(\dim(\Sigma)-2)$- face is either a circle or a line with two or three vertices and $\tilde\chi(\Gamma)=(-1)^{\dim\Gamma}$.}
\end{array}\right.
\end{array}$
\end{varprop}

\noindent
The above algebraic concepts of $\mbox{Bbm}_{_{{\bf k}}},$ $\mbox{CM}_{_{{\bf k}}}$ and $2\mbox{-CM}_{_{{\bf k}}}$-ness are strongly related to those for weak manifolds. Indeed, for triangulations of polytopes they are equivalent to
$n\mbox{-whm}_{_{{\bf k}}},  n\mbox{-jwhm}_{_{{\bf k}}}$ and
$n\mbox{-whsp}_{_{{\bf k}}}$ resp., implying that the following definitions are consistent with the original combinatorial/algebraic concepts:

\smallskip\noindent
{\bf Definition 7.} $X \mbox{ is Bbm}_{_{{\bf G}}}\ ({\rm CM}_{_{{\bf G}}},\ 2\mbox{-CM}_{_{{\bf G}}}) \mbox{ if } X \mbox{ is an } n\mbox{-whm}_{_{{\bf G}}}  (n\mbox{-jwhm}_{_{{\bf G}}}$,
$n\mbox{-whsp}_{_{{\bf G}}})$.\\
A simplicial complex $\Sigma$ is defined to be
${\rm Bbm}_{_{{\bf G}}}$, ${\rm CM}_{_{{\bf G}}}$ resp.
$2\mbox{-CM}_{_{{\bf G}}}$ if $|\Sigma|$ is.

\smallskip
For Bbm and CM, the consistency is easily checked through the homology-characterizations given above, while \cite{Fo1} provide a proof that the above definition of 2-CM-ness i consistent with Baklawski's original definition in \cite{Ba}.
In particular:
$\Sigma\mbox{ is 2-CM}_{_{{\bf k}}}$ if and only if $\Sigma$ is $\mbox{CM}_{_{{\bf k}}}$ and $\hat{\mbox{H}}_{n-1}({\rm cost}_{\Sigma}{\delta};\mathbf{k})=0$, $\text{ for all }\delta\in\Sigma,$  cf. \cite{Sta} p. 94.

\begin{varprop}~{\bf 4.}\label{Prop:04}
The following conditions are equivalent:

{\bf a.} $\Sigma$ is ${\rm Bbm}_{_{{\bf k}}}$,

{\bf b.} {\rm(Schenzel)} $\Sigma$ is pure and $\mathrm{Lk}_{\Sigma}\delta  \text{ is }
\mbox{CM}_{_{{\bf k}}}  \text{ for all }  \emptyset_{o}\ne\delta\in\Sigma,$

{\bf c.} {\rm(Reisner)} $\Sigma \text{ is pure and }
\mbox{Lk}_{\Sigma}v \text{ is } \mbox{CM}_{_{{\bf k}}} \text{for all } v\in V_{\Sigma}.$
\end{varprop}

{\bf Example.} When limited to compact polytopes  and a field {\bf k} as coefficient module/base ring, \cite{Sta} p. 73 gives the following Buchsbaum-equivalence using ring-theoretical {\it local cohomology};

{\bf d.} (Schenzel) $\Delta \text{ is Buchsbaum if and only if }
\dim{\mathrm{H}}^{i}_{{{\bf k}\langle\Sigma\rangle}_{{\bf +}}}
\!(\mathbf{k}\langle\Sigma\rangle)\le\infty\text{ if }0\le i < \dim\mathbf{k}\langle\Sigma\rangle),
$
in which case ${\mathrm{H}}^{i}_{{\bf k[\Sigma]}_{+}}({\mathbf{k}}\langle\Sigma\rangle)
\cong
{\hat{\mathbf{H}}}_{i-1}(|\Sigma|;\mathbf{k});$
for proof cf. \cite{S&V} p. 144.\\
Here, ``$\dim$" is the {\it Krull dimension}, which for Stanley-Reisner rings is simply ``1 + the simplicial dimension".

$\mbox{Given finite } \Sigma_{1},\ \Sigma_{2}\
 {\it n}\mbox{-jwhm}_{_{{\bf k}}} \mbox{ (equivalently, finite and CM}_{{\bf k}})$.
$\mbox{Since } \mathbf{k}
\langle\Sigma_{1}\rangle\otimes\mathbf{k}
\langle\Sigma_{2}\rangle
\cong
\mathbf{k}\langle\Sigma_{1}\ast \Sigma_{2}\rangle$
for Stanley-Reisner rings, the following {\it Künneth formula} for ring theoretical local cohomology is a direct consequence of the {\it Künneth formula for joins}; (``$\cdot_+$" indicates the unique homogeneous maximal ideal of ``$\cdot$", i.e., the irrelevant ideal of ``$\cdot$".)

\smallskip\noindent
{\bf 1.} ${\underline{\mathrm{H}}}^{q}_{_{
({\bf k}\langle\Sigma_{1}\rangle\otimes\mathbf{k}\langle\Sigma_{2}\rangle)_+}}
(\mathbf{k}{\langle\Sigma_{_{^{1}}}\rangle}\otimes\mathbf{k}{\langle\Sigma_{_{^{2}}}\rangle})
\cong
$
$
\bigoplus_{i+j=q}
\big({\underline{\mathrm{H}}}^{i}_{_{{\bf k}\langle\Sigma_{1}\rangle_+}}\!\!
(\mathbf{k}{\langle\Sigma _{_{^{1}}}\rangle})
\otimes
{\underline{\mathrm{H}}}^{j}_{_{{\bf k}\langle\Sigma_{_{^{2}}}\rangle_+}}\!\!
(\mathbf{k}{\langle\Sigma_{_{^{2}}}\rangle}) \big).
$

\medskip\noindent
{\bf 2.} Put:
$ \beta_{_{{\bf G}}}(X):= \mbox{inf}\{j \vert \exists x; x\in X \land
{\hat{\mathbf{H}}}_{j}(X,X\setminus_{o} x;\mathbf{G})\ne0\},$
where ``inf":="infimum".

\medskip
$\text{ For a finite }\Sigma, \mbox{  } \beta_{_{{\bf k}}}(\Sigma)$ is related to
the concepts ``depth of the ring $ \mathbf{k}\langle \Sigma \rangle $"
and ``C-M-ness of $\mathbf{k} \langle \Sigma \rangle $" through
$\beta_{_{{\bf k}}}(\Sigma)=\mbox{ depth } ( \mathbf{k}
\langle\Sigma\rangle)-1,$ in \cite{Br} Ex. 5.1.23 p. 214 and \cite{Sta} p. 142 Ex. 34. See also \cite{Mu}.

\subsubsection{Product and joins of Gorenstein complexes}Let ${\mathbf{k}}$ denote either a field or the integers ${\mathbb{Z}}.$
Proofs of the following two propositions are found in \cite{Fo2}.

\begin{varprop}~{\bf 5.}\label{Prop:05}
$\Sigma_{1}\ast\Sigma_{2}$ Gorenstein over ${\mathbf{k}}\Leftrightarrow
\Sigma_{1},\Sigma_{2}$ both Gorenstein over ${\mathbf{k}}$.
\end{varprop}

The (ordered simplicial) product of simplicial complexes is defined at the end of the chapter on Stanley-Reisner rings as Definition~\ref{defOrderSimpProd} p.~\pageref{defOrderSimpProd}.

\begin{varprop}~{\bf 6.}\label{Prop:06}
Let $\Sigma_1$ and $\Sigma_2$ be
two arbitrary simplicial complexes with $\dim\Sigma_i\!\geq1, \ (i=1,2)$ and  with linear orders on their vertex sets $V_{\Sigma_1}, V_{\Sigma_2}$ respectively, then,
$$\hskip-6.0cm\mbox{\rm (I) }\hskip3cm\Sigma_1\times \Sigma_2
\mbox{ Gorenstein over } {\mathbf{k}} $$
is equivalent to the disjunction of the following two statements;
$$\hskip-1.3cm\mbox{\rm (II) }{\Big\{}
\begin{array}{ll}\Sigma_i
\mbox{ are both Gorenstein over } {\mathbf{k}}
\mbox{ with exactly one cone point } v_i, i=1,2, \\
\mbox{ which either both are minimal } \mathbf{or}
\text{ both are maximal}.\end{array}$$
$$\mbox{\rm (III) }{\Big\{}\begin{array}{ll}\Sigma_i \text{ are both Gorenstein over }
{\mathbf{k}} \text{ with exactly two cone points } v_{ij}, 1\le i,  j\le 2, \\
\mbox{ where } v_{i1} \text{ are minimal elements in } V_{\Sigma_i}, \mbox{ and } v_{i2} \text{ are maximal in } V_{\Sigma_i}.\end{array}$$
\centerline{$\Bigl({\rm i.e.\ (I)} \Longleftrightarrow \mbox{\rm (II)} \lor \mbox{\rm (III)}  \Bigr)$}
\end{varprop}
{\bf Note.} The cases when either of the complexes have dimension less then 1 can be sorted out using Prop. 3. (iv) and (v) above.

\noindent\noindent\hrulefill
\vskip-0.4cm
\noindent\noindent\hrulefill
\vskip-0.15cm

\smallskip

\underbar{External} hyperlinks in \underbar{pdf-files} should be right-clicked to decide which action is preferred.

\vskip-0.3cm
\noindent\noindent\hrulefill
\vskip-0.4cm
\noindent\noindent\hrulefill
\vskip-0.15cm

\end{document}